\documentclass[a4paper,leqno,11pt]{article}
\usepackage{amssymb,amsfonts,amsmath,amsthm}
\usepackage{epsfig,epsf,curves}
\usepackage{bbm,epic,eepic}
\usepackage[mathscr]{eucal}
\theoremstyle{remark}

\begin{document}

\sloppy


\def\Mat#1#2#3#4{\left(\!\!\!\begin{array}{cc}{#1}&{#2}\\{#3}&{#4}\\ \end{array}\!\!\!\right)}

\def\tAA{{\Bbb A}}
\def\CC{{\Bbb C}}
\def\HH{{\Bbb H}}
\def\NN{{\Bbb N}}
\def\QQ{{\Bbb Q}}
\def\RR{{\Bbb R}}
\def\ZZ{{\Bbb Z}}
\def\FF{{\Bbb F}}
\def\SS{{\Bbb S}}
\def\GG{{\Bbb G}}
\def\PP{{\Bbb P}}
\def\LL{{\Bbb L}}

\title{On Truncation of irreducible representations of Chevalley groups}

\author{Joachim Mahnkopf}

\maketitle

{\bf Abstract } We prove part of a higher rank analogue of the Gouvea-Mazur Conjecture (cf. [G-M]). More precisely, let $\tilde{\bf G}$ be a connected, 
reductive ${\Bbb Q}$-split group and $\Gamma$ an arithmetic subgroup of $\tilde{\bf G}$. We show that the dimension of the slope $\alpha$ subspace of the 
cohomology of $\Gamma$ with values in an irreducible $\tilde{\bf G}$-module $L$ is bounded independently of $L$. The proof is based on general principles 
of the representation theory of algebraic groups; in particular, we study truncations of highest weight modules of Chevalley groups.


\vspace{0.8cm}

\centerline{\bf \large Introduction}

{\bf (0.1) } The first boundedness result for the dimension of the slope subspaces of the cohomology groups of arithmetic subgroups has been obtained by Hida; he showed that the dimension of the slope $0$-subspace of the cohomology of an arithmetic subgroup $\Gamma$ in ${\rm GL}_2({\Bbb Z})$ with 
coefficients in the irreducible representation $L_k$ of highest weight $k\in{\Bbb N}_0$ is bounded independently of $k$ (cf. e.g. [Hi], chapter 7.2). 
Hida even showed that this dimension is constant as a function of $k$ and later he generalized his results to higher rank. Following a suggestion by 
R. Taylor, Buzzard extended Hida's result to spaces of arbitrary slope (but still considering arithmetic subgroups of ${\rm GL}_2$) (cf. [Bu]). A. Pande, following the method of Buzzard/Taylor, further extended Buzzard's result to ${\rm GL}_2$ over imaginary quadratic fields (cf. [P]). 
We also mention the recent work by Harder in the slope $0$ case (cf. [H]).

\medskip

{\bf (0.2) } In this article we prove boundedness of the dimension of the slope subspaces for arithmetic subgroups in arbitrary reductive algebraic groups. 
To describe this in more detail, let $\tilde{\bf G}$ be a connected, reductive, ${\Bbb Q}$-split group with derived (semi simple) group ${\bf G}$. We 
choose a maximal ${\Bbb Q}$-split torus $\tilde{\bf T}$ in $\tilde{\bf G}$ and a maximal ${\Bbb Q}$-split torus ${\bf T}\le\tilde{\bf T}$ in ${\bf G}$ 
and we denote by $L_{\tilde{\lambda}}$ the irreducible algebraic $\tilde{\bf G}$-module of highest weight $\tilde{\lambda}\in X(\tilde{\bf T})$. We note 
that $L_{\tilde{\lambda}}$ is defined over ${\Bbb Z}$. We fix a prime $p$ 
and we let $\Gamma\le {\bf G}({\Bbb Z})$ be an arithmetic subgroup. As in the above mentioned works, we assume that $\Gamma$ satisfies a certain level 
condition at $p$ (cf. section (4.2) for a precise definition). The group $\Gamma$ acts on $L_{\tilde{\lambda}}({\Bbb Q}_p)$ and on 
$L_{\tilde{\lambda}}({\Bbb Z}_p)$ and we define the slope $\alpha$ 
subspace of the group cohomology $H^i(\Gamma,L_{\tilde{\lambda}}({\Bbb Q}_p))$ with respect to a normalized Hecke operator at $p$. Our main result then is as 
follows (cf. section (6.6)). 

\medskip

{\bf Theorem. }  {\it Let $\beta\in{\Bbb Q}_{\ge 0}$ and select $i\in{\Bbb N}_0$. There is a constant $C=C(\beta,i)$ such that 
$$
\sum_{0\le\alpha\le\beta} {\rm dim}\, H^i(\Gamma,L_{\tilde{\lambda}}({\Bbb Q}_p))^\alpha\le C
$$
for all dominant weights $\tilde{\lambda}$.
}

\medskip

{\bf (0.3) } The proof of the above theorem is based on consideration of truncations of irreducible representations of ${\bf G}$. 
To explain this we let $\tilde{\lambda}\in X(\tilde{\bf T})$ be a dominant algebraic weight and we denote by $\lambda^\circ=\tilde{\lambda}|_{\bf T}$ resp. 
$\lambda={\rm d}\,\lambda^\circ \in{\rm Lie}({\bf T})^*$ the restriction of $\tilde{\lambda}$ to ${\bf T}$ resp. the differential of $\lambda^\circ$.  We 
denote by $L_{\lambda}$ the irreducible ${\rm Lie}({\bf G})$-module of highest weight $\lambda$, hence, ${\bf G}$ acts on $L_{\lambda}$ and $L_\lambda$ 
is the algebraic irreducible ${\bf G}$-module of highest weight $\lambda^\circ$ and is defined over ${\Bbb Z}$. 
%
%
The truncation ${\bf L}^{[r]}_\lambda({\Bbb Z}_p)$ of length $r$ of the $\Gamma$-module $L_{{\lambda}}$ is defined as a quotient of 
${\Bbb Z}_p$-modules $L_{{\lambda}}({\Bbb Z}_p)/L_{{\lambda}}({\Bbb Z}_p,r)$, where $L_{{\lambda}}({\Bbb Z}_p,r)\le L_{{\lambda}}({\Bbb Z}_p)$ 
is a certain ${\Bbb Z}_p$-submodule, which in particular contains all weight spaces whose weight has 
"relative height" strictly larger than $r$ (cf. (1.3)). Thus, we truncate from $L_{{\lambda}}({\Bbb Z}_p)$ all "sufficiently high" weight spaces. 
${\bf L}^{[r]}_{{\lambda}}({\Bbb Z}_p)$ no longer is a ${\bf G}({\Bbb Z}_p)$-module (or ${\bf G}({\Bbb Z})$-module), but it still is a $\Gamma$-module 
of finite cardinality. 
The Theorem then is a consequence of the following two properties of the truncation.

\medskip

{\bf Proposition A } (cf. (6.5) Proposition). {\it Let $r$ be an integer bigger than $\alpha+1$. Then 
$$
{\rm dim}\, H^i(\Gamma,L_{\tilde{\lambda}}({\Bbb Q}_p))^\alpha\le v_p\big(\#(H^i(\Gamma,{\bf L}_\lambda^{[r]}({\Bbb Z}_p)))\big).
$$ 
}
This follows from a closer analysis of the Hecke Operator acting on the cohomology of $\Gamma$ in section 5 
together with a cohomological argument due to Hida (cf. the proof of (6.4) Proposition or [Hi] chapter 7.2).

%
%

\medskip

{\bf Proposition B } (cf. (4.4) Proposition). {\it  The isomorphism class of the $\Gamma$-module ${\bf L}^{[r]}_\lambda({\Bbb Z}_p)$ 
only depends on ${\lambda}$ modulo $p^{\lceil\frac{p}{p-1} \rceil r}\Gamma_{\rm sc}$ ($\Gamma_{\rm sc}=$ the weight lattice of ${\rm Lie}({\bf G})$).

}


\medskip

In the above form Proposition B only holds if the weight $\lambda$ is sufficiently regular; if $\lambda$ comes close to the boundary of the 
Weyl chamber the statement has to be modified (cf. Proposition (4.4) for a precise statement). We note that in the main part we will consider somewhat more 
generally truncations of $L_{{\lambda}}({\Bbb Z})$.

\medskip

We deduce Proposition B from its differential version for irreducible representations of semi simple Lie algebras by using the method of Chevalley groups.
We note that truncations of irreducible ${\rm GL}_2$-modules appear in the above mentioned works of Buzzard and Pande as certain subspaces of symmetric 
powers of the standard representation of ${\rm GL}_2$. Our definition is an extension to the higher rank case based on the semi simple representation theory of 
Chevalley groups. 
Also, in the context of symmetric powers of the standard representation of ${\rm GL}_2$ a result 
analogous to the above Proposition appears in [P].

\medskip

{\bf (0.4) } 
Buzzard and Pande make use of 
the fact that in the ${\rm GL}_2$-case the (interesting) cohomology appears in degree $1$, hence, any cocycle 
already is determined by its values on a set of generators of $\Gamma$. 
Based on this they even obtain an explicit upper bound for the dimension of the slope spaces and in 
the case of the trivial congruence subgroup $\Gamma$ in a indefinite quaternion algebra Pande also proves local constancy of the dimension of the 
slope spaces. We hope to be able to deal with these questions in the context of this article in the future.



\medskip

{\bf (0.5) } Our motivation for writing this article comes from the Gouvea-Mazur Conjecture. On the one hand, the above Theorem confirms a part of a higher 
rank analogue of the Gouvea-Mazur Conjecture. On the other hand, we hope that it will be an ingredient in an approach to the construction of $p$-adic 
families of modular forms (i.e. the second part of the Gouvea-Mazur Conjecture), which is based on a more elementary comparison of trace formulas 
and does not make use of advanced theories like rigid analytic geometry (cf. [M 1,2]).

\bigskip

{\bf Acknowledgement. } I am grateful to the referee for pointing out errors in earlier versions and for several helpful remarks and suggestions.

\bigskip

\section{Truncation of irreducible representations of semi simple Lie algebras}

{\bf (1.1) Notations. } We fix a complex semi-simple Lie algebra ${\mathfrak g}$ of rank $\ell$. We introduce the following notations.

\begin{itemize}

\item ${\mathfrak h}$ denotes a fixed a Cartan subalgebra of ${\mathfrak g}$ and $\Phi$ is the set of roots of ${\mathfrak g}$ with respect to ${\mathfrak h}$.

\item $\Delta$ is a fixed choice of a basis of $\Phi$ and $\Phi^+$ resp. $\Phi^-$ is the set of positive resp. negative roots in $\Phi$. 

\item ${\mathfrak g}(\alpha)$ is the weight $\alpha$ subspace of ${\mathfrak g}$ (with respect to the adjoint action of ${\mathfrak h}$) and we set ${\mathfrak n}=\bigoplus_{\alpha\in \Phi^+}{\mathfrak g}(\alpha)$, ${\mathfrak n}^-=\bigoplus_{\alpha\in \Phi^-}{\mathfrak g}(\alpha)$.

\item $(\quad,\quad)$ is the Killing form on ${\mathfrak g}$ and for $\lambda\in {\mathfrak h}^*$ we define $t_\lambda\in{\mathfrak h}$ by
$(t_\lambda,h)=\lambda(h)$ for all $h\in {\mathfrak h}$. 
In this way, $(\quad,\quad)$ induces a pairing on ${\mathfrak h}^*$ and for
$\lambda,\mu\in{\mathfrak h}^*$ we set $\langle \lambda,\mu \rangle=2(\lambda,\mu)/(\mu,\mu)$.

\item for every $\alpha\in \Phi$ we define the coroot $h_\alpha=\frac{2 t_\alpha}{(t_\alpha,t_\alpha)}\in{\mathfrak h}$; hence, $\lambda(h_\alpha)=\langle \lambda,\alpha \rangle$. Moreover, we denote by $\{\omega_\alpha\}_{\alpha\in \Delta}$ the set of fundamental weights, i.e. $\omega_\alpha(h_\beta)=\delta_{\alpha,\beta}$ for all $\alpha,\beta\in \Delta$.


\item we denote by $\Gamma_{\rm sc}$ the weight lattice of ${\mathfrak g}$ consisting of all integral weights in 
${\mathfrak h}^*$, i.e. $\Gamma_{\rm sc}$ consists of all weights
$\lambda\in{\mathfrak h}^*$ such that $\lambda(h_\alpha)\in{\Bbb Z}$ for all $\alpha\in\Phi$. $\Gamma_{\rm ad}$ is the root lattice, i.e. the subgroup 
of ${\mathfrak h}^*$ generated by the roots. We note that $\Gamma_{\rm sc}\le {\mathfrak h}^*$ is a ${\Bbb Z}$-lattice in ${\mathfrak h}^*$ with basis $\{\omega_\alpha\}_{\alpha\in\Delta}$ and $\Gamma_{\rm ad}\le \Gamma_{\rm sc}$ is a sublattice with basis $\Delta$. We write
$\lambda\ge \mu$ ($\lambda,\mu\in\Gamma_{\rm sc}$) if $\lambda-\mu$ is a linear combination of elements in $\Phi^+$ (or, equivalently, in $\Delta$) with positive coefficients (cf. [Hu], p. 47). 
%
We denote by 
$$
{\rm ht}={\rm ht}_\Delta:\,\Gamma_{\rm ad}\rightarrow {\Bbb Z}
$$ 
the height function on $\Gamma_{\rm ad}$, i.e. ${\rm ht}(\lambda)=\sum_{\alpha\in \Delta} c_\alpha$ for $\lambda=\sum_{\alpha\in\Delta} c_\alpha \alpha\in\Gamma_{\rm ad}$.


\item ${\cal U}$ is the universal enveloping algebra of ${\mathfrak g}$ and ${\cal U}^+$, resp. ${\cal U}^-$ resp. ${\cal U}^o$ is the
universal enveloping algebra of ${\mathfrak n}$ resp. ${\mathfrak n}^-$ resp. ${\mathfrak h}$ .

\item For any $\alpha\in\Phi$ we choose a root vector $x_\alpha\in {\mathfrak g}(\alpha)$ such that $\{x_\alpha,\,\alpha\in
\Phi,\;h_\alpha,\,\alpha\in \Delta\}$ is a {\it Chevalley basis} of ${\mathfrak g}$. 

\item We fix an ordering $\{\alpha_1,\ldots,\alpha_s\}$ of $\Phi^+$; for any multi index 
${\bf n}=(n_1,\ldots,n_s)\in{\Bbb N}_0^s$ we set
$$
X_+^{\bf n}=\frac{x_{\alpha_1}^{n_1}}{n_1!}\cdots \frac{x_{\alpha_s}^{n_s}}{n_s!}\in{\cal U}^+
$$
and 
$$
X_-^{\bf n}=\frac{x_{-\alpha_1}^{n_1}}{n_1!}\cdots \frac{x_{-\alpha_s}^{n_s}}{n_s!}\in{\cal U}^-.
$$
Moreover, for ${\bf n}=(n_\alpha)_{\alpha\in\Delta}\in{\Bbb N}_0^\Delta$ we set
$$
H^{\bf n}=\prod_{\alpha\in \Delta}{h_\alpha\choose n_\alpha}\in{\cal U}^o,
$$
where ${h\choose n}=h(h-1)\cdots (h-n+1)/n!$ (${h\choose 0}=1$). 

\item We denote by ${\cal U}_{\Bbb Z}\le {\cal U}$ the ${\Bbb Z}$-algebra generated by the elements $x_\alpha^n/n!$ ($\alpha\in\Phi$, $n\in{\Bbb N}_0$). Theorem 26.4 in [Hu], p. 156 implies that ${\cal U}_{\Bbb Z}$ is a lattice in ${\cal U}$ with basis 
$$
{\cal U}_{\Bbb Z}=\bigoplus_{{\bf n}_1\in{\Bbb N}_0^{s},{\bf n}_2\in{\Bbb N}_0^{\Delta},{\bf n}_3\in{\Bbb N}_0^{s}} 
{\Bbb Z}\, X_-^{{\bf n}_1}H^{{\bf n}_2}X_+^{{\bf n}_3}.
$$
Thus, ${\cal U}_{\Bbb Z}$ is a ${\Bbb Z}$-form of the associative algebra ${\cal U}$. Similarly, ${\cal U}^+$ resp. ${\cal U}^-$
resp. ${\cal U}^o$ have ${\Bbb Z}$-forms ${\cal U}^+_{\Bbb Z}$ resp. ${\cal U}^-_{\Bbb Z}$
resp. ${\cal U}^o_{\Bbb Z}$ with ${\Bbb Z}$-basis $\{X_+^{{\bf n}},\; {\bf n}\in{\Bbb N}_0^{s}\}$ resp. 
$\{X_-^{\bf n},\; {\bf n}\in{\Bbb N}_0^{s}\}$ resp. $\{H^{\bf n},\; {\bf n}\in{\Bbb N}_0^{\Delta}\}$. 
In particular, we obtain ${\cal U}_{\Bbb Z}={\cal U}_{\Bbb Z}^-{\cal U}_{\Bbb Z}^o{\cal U}_{\Bbb Z}^+$.

\end{itemize}

We denote by $\lceil x\rceil$ the smallest integer, which is equal to or bigger than $x$. We define the binomial coefficient ${z\choose k}=\frac{z(z-1)\cdots(z-k+1)}{k!}$, $z\in{\Bbb Z}$, $k\in{\Bbb N}$ and we put ${z\choose 0}=1$. It is easy to see that ${z\choose k}\in{\Bbb Z}$ for all $z\in{\Bbb Z}$ 
and $k\in{\Bbb N}_0$. 
We denote by $v_p$
the $p$-adic valuation on $\bar{\Bbb Q}_p$ normalized by $v_p(p)=1$. Moreover, we fix an embedding $\bar{\Bbb Q}\hookrightarrow \bar{\Bbb Q}_p$, hence, $v_p$ defines a $p$-adic valuation on $\bar{\Bbb Q}$.

\bigskip

{\bf (1.2) Irreducible representations of semi simple Lie algebras. } Let $\lambda\in \Gamma_{\rm sc}$ be an integral and dominant weight. We denote by $(\rho_{\lambda},L_{\lambda})$ the finite dimensional, irreducible complex representation of 
${\mathfrak g}$ of highest weight $\lambda$. We note that $L_\lambda$ is an ${\cal U}$-module. We denote by $L_{\lambda}(\mu)$ the weight $\mu$ subspace of 
$L_{\lambda}$, $\Pi_\lambda$ is the set of all nontrivial weights $\mu$ of $L_\lambda$ (i.e $L_\lambda(\mu)\not=0$) and $\Gamma_{\lambda}\le{\mathfrak h}^*$ is 
the weight lattice of $\rho_\lambda$, i.e. $\Gamma_\lambda$ is the subgroup of $\Gamma_{\rm sc}$ generated by $\Pi_\lambda$. Thus, if $L_\lambda$ is faithful we have
$$
\Gamma_{\rm sc}\ge \Gamma_\lambda\ge \Gamma_{\rm ad}.
$$
(cf. [B], ch. 2.1, equation (2), p. 6). 
We fix a maximal vector $v_\lambda\in L_\lambda$, i.e. $v_\lambda$ has weight $\lambda$ and $L_\lambda={\cal U}^- v_\lambda$, and we set $L_\lambda({\Bbb Z})={\cal U}^-_{\Bbb Z}v_\lambda$. Theorem 27.1 b.) in [Hu], p. 158 and its proof imply that $L_\lambda({\Bbb Z})$ is an admissible, i.e. ${\cal U}_{\Bbb Z}$-invariant 
lattice in $L_\lambda$;
%
%
moreover, Theorem 27.1 a.) in [Hu], p. 158 implies 
$$
L_\lambda({\Bbb Z})=\bigoplus_{\mu\in \Pi_\lambda} L_\lambda({\Bbb Z},\mu),
$$
where $L_\lambda({\Bbb Z},\mu)=L_\lambda(\mu)\cap L_\lambda({\Bbb Z})$ is the weight $\mu$ subspace in $L_\lambda({\Bbb Z})$.
For any ${\Bbb Z}$-algebra $R$ we set $L_\lambda(R)=L_\lambda({\Bbb Z})\otimes_{\Bbb Z}R$ and $L_\lambda(R,\mu)
=L_\lambda({\Bbb Z},\mu)\otimes_{\Bbb Z}R$. Hence,
$$
L_\lambda(R)=\bigoplus_{\mu\in \Pi_\lambda} L_\lambda(R,\mu).
$$

\bigskip

{\bf (1.3) The truncating submodule of an irreducible representation. } In the remainder of section 1 we fix a prime $p\in{\Bbb N}$ and we define the following ("Iwahori-type") ${\Bbb Z}$-subalgebra of ${\cal U}_{\Bbb Z}$:
$$
{\cal S}={\Bbb Z}[\frac{x_{-\alpha}^n}{n!},\,\alpha\in\Phi^+,n\in{\Bbb N}_0,\,p^{m {\rm ht}(\alpha)}\frac{x_\alpha^m}{m!},\,\alpha\in\Phi^+,m\in{\Bbb N}_0]\,\le{\cal U}_{\Bbb Z},
$$
i.e. ${\cal S}$ is the ${\Bbb Z}$-subalgebra of ${\cal U}_{\Bbb Z}$ generated by the elements $\frac{x_{-\alpha}^n}{n!}$, $\alpha\in\Phi^+$, $n\in{\Bbb N}_0$
 and $p^{m {\rm ht}(\alpha)}\frac{x_\alpha^m}{m!}$, $\alpha\in\Phi^+$, $m\in{\Bbb N}_0$. We select an integral and dominant weight 
$\lambda\in \Gamma_{\rm sc}$. $L_\lambda({\Bbb Z})$ then becomes an ${\cal S}$-module via restriction. In the following we want to introduce 
the truncation of the ${\cal S}$-module $L_\lambda({\Bbb Z})$. We define a (relative) height function  
$$
{\rm ht}_\lambda:\,\Pi_\lambda\rightarrow{\Bbb N}_0
$$ 
on the set of weights of $L_\lambda$ by 
$$
{\rm ht}_\lambda(\mu)={\rm ht}(\lambda-\mu)
$$ 
(note that $\lambda-\mu\in\Gamma_{\rm ad}$, hence, ${\rm ht}(\lambda-\mu)$ is defined). Explicitly, if $\mu=\lambda-\sum_{\alpha\in\Delta} c_\alpha \alpha$, $c_\alpha\in{\Bbb N}_0$ (note that any $\mu\in\Pi_\lambda$ has this form), then 
${\rm ht}_\lambda(\mu)=\sum_{\alpha\in\Delta} c_\alpha$. 
%
%
We also note that
$$
{\rm ht}_\lambda(\mu+\omega)={\rm ht}_\lambda(\mu)-{\rm ht}(\omega)\leqno(1)
$$
for all $\mu\in \Pi_\lambda$ and all $\omega\in \Gamma_{\rm ad}$ such that $\mu+\omega$ again is contained in $\Pi_\lambda$. 
For any integer $r\in{\Bbb N}_0$ we then define the following ${\Bbb Z}$-submodule of the highest weight module $L_\lambda({\Bbb Z})$: 
$$
L_\lambda({\Bbb Z},r)=\bigoplus_{\mu\in\Pi_\lambda \atop 0\le {\rm ht}_\lambda(\mu)\le r} p^{r-{\rm ht}_\lambda(\mu)}\,L_\lambda({\Bbb Z},\mu) \oplus
\bigoplus_{\mu\in\Pi_\lambda \atop {\rm ht}_\lambda(\mu)> r} L_\lambda({\Bbb Z},\mu).
$$
Again, we put $L_\lambda(R,r)=L_\lambda({\Bbb Z},r)\otimes_{\Bbb Z}R$ for any ${\Bbb Z}$-algebra $R$, hence, 
$$
L_\lambda(R,r)=\bigoplus_{\mu\in\Pi_\lambda \atop 0\le {\rm ht}_\lambda(\mu)\le r} p^{r-{\rm ht}_\lambda(\mu)}\,L_\lambda(R,\mu) \oplus
\bigoplus_{\mu\in\Pi_\lambda \atop {\rm ht}_\lambda(\mu)> r} L_\lambda(R,\mu).
$$
We note that $L_\lambda({\Bbb Z},r)$ is a ${\Bbb Z}$-lattice in $L_\lambda({\Bbb Q},r)$.

\bigskip

{\bf (1.4) Lemma. }{\it $L_\lambda({\Bbb Z},r)$ is an ${\cal S}$-invariant submodule of $L_\lambda({\Bbb Z})$.

}

\medskip

{\it Proof. } Since $\frac{x_{-\alpha}^n} {n!} L_\lambda({\Bbb Z})\subseteq L_\lambda({\Bbb Z})$ we obtain for all $\mu\in\Pi_\lambda$ and $\alpha\in \Phi^+$ 
$$
{\frac{x_{-\alpha}^n}{n!}}\, L_\lambda({\Bbb Z},\mu)\subseteq L_\lambda(\mu-n\alpha)\cap L_\lambda({\Bbb Z})
=L_\lambda({\Bbb Z},\mu-n\alpha).
$$
Since ${\rm ht}_\lambda(\mu-n\alpha)\ge {\rm ht}_\lambda(\mu)$ (cf. equation (1)) the definition of $L_\lambda({\Bbb Z},r)$
immediately implies that ${\frac{x_{-\alpha}^n}{n!}}\, L_\lambda({\Bbb Z},r)\subseteq L_\lambda({\Bbb Z},r)$.

\medskip

To show that all generators $p^{n {\rm ht}(\alpha)}\frac{x_\alpha^n}{n!}$ with $\alpha\in\Phi^+$ leave 
$L_\lambda({\Bbb Z},r)$ invariant, we distinguish cases.  

1. First, we consider weights $\mu\in\Pi_\lambda$ satisfying ${\rm ht}_\lambda(\mu)\le r$. We let $v\in p^{r-{\rm ht}_\lambda(\mu)}\,L_\lambda({\Bbb Z},\mu)$ be arbitrary. For all
$\alpha\in \Phi^+$ we obtain using equation (1)
$$
p^{n {\rm ht}(\alpha)} \frac{x_\alpha^n}{n!}(v)\in  p^{n {\rm ht}(\alpha)} p^{r-{\rm ht}_\lambda(\mu)} L_\lambda({\Bbb Z},\mu+n\alpha)
=p^{r-{\rm ht}_\lambda(\mu+n\alpha)}\,L_\lambda({\Bbb Z},\mu+n\alpha).
$$
Since ${\rm ht}_\lambda(\mu+n\alpha)\le {\rm ht}_\lambda(\mu)\le r$ we deduce that $p^{r-{\rm ht}_\lambda(\mu+n\alpha)} L_\lambda({\Bbb Z},\mu+n\alpha)$ is contained in 
$L_\lambda({\Bbb Z},r)$. Since $v$ was arbitrary this implies that $p^{n {\rm ht}(\alpha)} \frac{x_\alpha^n}{n!} p^{r-{\rm ht}_\lambda(\mu)} L_\lambda({\Bbb Z},\mu)\subseteq L_\lambda({\Bbb Z},r)$ (note that we may assume that $\mu+n\alpha\in\Pi_\lambda$ since otherwise $x_\alpha^n v=0$).  

2. Second, we consider weights $\mu\in\Pi_\lambda$ with ${\rm ht}_\lambda(\mu)>r$. We let $v\in L_\lambda({\Bbb Z},\mu)$ be arbitrary. As above we find
$$
p^{n {\rm ht}(\alpha)}  \frac{x_\alpha^n}{n!} (v)\in p^{n{\rm ht}(\alpha)} \,L_\lambda({\Bbb Z},\mu+n\alpha).\leqno(2)
$$
If ${\rm ht}_\lambda(\mu+n\alpha)> r$ then $p^{n{\rm ht}(\alpha)} \,L_\lambda({\Bbb Z},\mu+n\alpha)$ 
clearly is contained in $L_\lambda({\Bbb Z},r)$. Thus, we may assume that 
${\rm ht}_\lambda(\mu+n\alpha)\le r$. Since ${\rm ht}_\lambda(\mu)> r$ we obtain using equation (1)
$$
r-{\rm ht}_\lambda(\mu+n\alpha)=r-{\rm ht}_\lambda(\mu)+n{\rm ht}(\alpha)\le n{\rm ht}(\alpha).
$$
Thus, the definition of $L_\lambda({\Bbb Z},r)$ shows that $p^{n{\rm ht}(\alpha)} \,L_\lambda({\Bbb Z},\mu+n\alpha)$
is contained in $L_\lambda({\Bbb Z},r)$. Hence, in the second case we also obtain that
$p^{n {\rm ht}(\alpha)} \frac{x_\alpha^n}{n!} L_\lambda({\Bbb Z},\mu)\subseteq  L_\lambda({\Bbb Z},r)$. 

Thus, cases 1 and 2 show that $p^{n {\rm ht}(\alpha)} \frac{x_\alpha^n}{n!} L_\lambda({\Bbb Z},r)\subseteq  L_\lambda({\Bbb Z},r)$ and the proof of the Lemma is 
complete.

\bigskip

{\bf (1.5) Truncation of an irreducible representation. } Since $L_\lambda({\Bbb Z},r)$ is invariant under ${\cal S}$ we obtain a representation 
$$
\rho_{\lambda}:\,{\cal S}\rightarrow {\rm End}(L_\lambda({\Bbb Z})/L_\lambda({\Bbb Z},r)),
$$
which we call the truncation of length $r$ of the representation 
$(\rho_\lambda,L_\lambda)$. We also call the ${\cal S}$-module 
$$
{\bf L}_\lambda^{[r]}(R)=L_\lambda({R})/L_\lambda({R},r)
$$ 
the truncation of length $r$ of the highest weight module $L_\lambda({R})$ and $L_\lambda({R},r)$ the truncating submodule. In the following we 
will only need to consider truncations in the cases $R={\Bbb Z}$ and $R={\Bbb Z}_p$. We note that
$$
{\bf L}^{[r]}_\lambda({\Bbb Z})\cong\bigoplus_{\mu\in \Pi_\lambda \atop {\rm ht}_\lambda(\mu)\le r} 
\frac{L_\lambda({\Bbb Z},\mu)}{p^{r-{\rm ht}_\lambda(\mu)} L_\lambda({\Bbb Z},\mu)}.\leqno(3)
$$
In particular, ${\bf L}^{[r]}_\lambda({\Bbb Z})$ is a finitely generated torsion ${\Bbb Z}$-module, which is annihilated by $p^r$.

\newpage

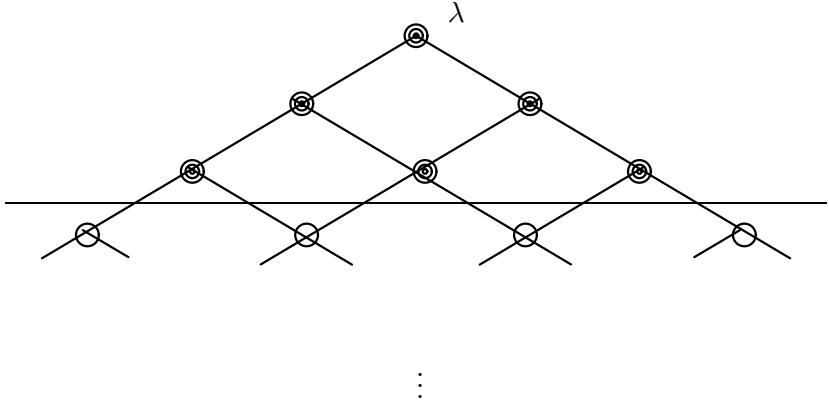
\begin{picture}(5,9)
\setlength{\unitlength}{0.6cm}

\thicklines
\put(12,0){\line(-5,-3){8.2}}
\put(12,0){\line(5,-3){8.2}}
\put(9.3,-1.4){\line(5,-3){6.1}}
\put(14.7,-1.4){\line(-5,-3){6.1}}
\put(17,-2.9){\line(-5,-3){3.6}}
\put(7,-2.9){\line(5,-3){3.6}}
\put(19.1,-4.3){\line(-5,-3){1}}
\put(4.7,-4.3){\line(5,-3){1}}

\put(12,0){\circle{0.5}}
\put(12,0){\circle{0.3}}
\put(12,0){\circle{0.1}}

\put(9.5,-1.5){\circle{0.5}}
\put(9.5,-1.5){\circle{0.3}}
\put(9.5,-1.5){\circle{0.1}}

\put(14.5,-1.5){\circle{0.5}}
\put(14.5,-1.5){\circle{0.3}}
\put(14.5,-1.5){\circle{0.1}}

\put(16.9,-3){\circle{0.5}}
\put(16.9,-3){\circle{0.3}}
\put(16.9,-3){\circle{0.1}}

\put(7.1,-3){\circle{0.5}}
\put(7.1,-3){\circle{0.3}}
\put(7.1,-3){\circle{0.1}}

\put(12.2,-3){\circle{0.5}}
\put(12.2,-3){\circle{0.3}}
\put(12.2,-3){\circle{0.1}}

\put(9.6,-4.4){\circle{0.5}}
\put(14.4,-4.4){\circle{0.5}}
\put(4.8,-4.4){\circle{0.5}}
\put(19.2,-4.4){\circle{0.5}}


\thinlines
\put(3,-3.7){\line(1,0){18}}


\put(12,-8){\vdots}

\put(2,-10.2){Root system $A_2$: $"\quad"=$ Weights of the truncation $L_\lambda({\Bbb Z})/L_\lambda({\Bbb Z},2)$ of length $2$}
\put(9.2,-11){in case $m_{\alpha_1}$, $m_{\alpha_2}>2$ \quad ($\Delta=\{\alpha_1,\alpha_2\}$).}

\put(7.6,-10){\circle{0.5}}
\put(7.6,-10){\circle{0.3}}
\put(7.6,-10){\circle{0.1}}


\put(12.7,0.3){$\lambda$}

\end{picture}


\vspace{8cm}

\section{Truncation of Verma modules}

{\bf (2.1) Integral Verma  modules. } In this section, we introduce an analogue of the truncation of irreducible representations for Verma modules. This will be used in section 3. We let $\lambda\in{\mathfrak h}^*$ and we denote by
$$
I_\lambda=\langle {\cal U}^+,h_\alpha-\lambda(h_\alpha),\,\alpha\in\Delta  \rangle_{\cal U}
$$
the left ideal (i.e. the ${\cal U}$-submodule) in ${\cal U}$ generated by ${\cal U}^+$ and $h_\alpha-\lambda(h_\alpha)$ for $\alpha\in\Delta$. The Verma 
module of highest weight $\lambda$ is defined as $V_\lambda:={\cal U}/I_\lambda$. Thus, $V_\lambda$ is a cyclic ${\cal U}$-module generated by the maximal 
vector ${\rm v}_\lambda=1+I_\lambda$. Using the PBW Theorem we even find that $V_\lambda={\cal U}^-{\rm v_\lambda}$ is a free ${\cal U}^-$-module, hence,
$$
V_\lambda=\bigoplus_{{\bf n}\in{\Bbb N}_0^s} {\Bbb C} X_-^{\bf n} {\rm v}_\lambda.\leqno(1)
$$
%
%
We denote by $V_\lambda(\mu)$ 
the weight $\mu$ subspace in $V_\lambda$. Hence, 
$$
V_\lambda=\bigoplus_{\mu\le \lambda} V_\lambda(\mu)\leqno(2)
$$ 
and
$$
V_\lambda(\mu)=\bigoplus_{{\bf n}\in{\Bbb N}_0^s\atop n_1\alpha_1+\cdots+n_s\alpha_s=\lambda-\mu} {\Bbb C} X_-^{\bf n} {\rm v}_\lambda\leqno(3)
$$ 
(to see this, note that by equation (1) $V_\lambda$ is the direct sum of the subspaces appearing on the right hand side of equation (3) and that these spaces have weight $\mu$). From now on, we assume in addition that $\lambda\in\Gamma_{\rm sc}$ is integral. We attach to $\lambda$ a ${\Bbb Z}$-linear mapping
$$
\lambda:\,{\cal U}_{\Bbb Z}^\circ=\bigoplus_{{\bf n}\in{\Bbb N}_0^\Delta} {\Bbb Z} H^{\bf n}\rightarrow{\Bbb Z} 
$$
by defining the values on basis elements as
$$
\lambda(H^{\bf n})=\prod_{\beta\in\Delta} {\lambda(h_\beta)\choose n_\beta}\qquad ({\bf n}=(n_\beta)_{\beta\in\Delta}\in{\Bbb N}_0^\Delta).
$$
We note that the integrality of $\lambda$ implies that $\lambda(h_\beta)=\langle\lambda,\beta\rangle\in{\Bbb Z}$ for all $\beta\in \Delta$, hence, $\lambda(H^{\bf n})\in{\Bbb Z}$. 
(Later we will in addition assume that $\lambda$ is dominant which implies that $\lambda(h_\beta)\in{\Bbb N}_0$). 

\medskip

We define the integral Verma module as 
$$
V_\lambda({\Bbb Z})={\cal U}_{\Bbb Z} {\rm v}_\lambda.
$$
%
%
A straightforward computation yields
$$
H{\rm v}_\lambda=\lambda(H){\rm v}_\lambda
$$ 
for all $H\in{\cal U}_{\Bbb Z}^\circ$. Since ${\cal U}_{\Bbb Z}^+$ annihilates ${\rm v}_\lambda$ we thus obtain
$$
V_\lambda({\Bbb Z})={\cal U}_{\Bbb Z}^- {\rm v}_\lambda=\bigoplus_{{\bf n}\in{\Bbb N}_0^s} {\Bbb Z} X_-^{\bf n} {\rm v}_\lambda;\leqno(4)
$$
hence, $V_\lambda({\Bbb Z})$ is a free ${\cal U}_{\Bbb Z}^-$-module and a ${\Bbb Z}$-lattice in $V_\lambda$. 
Moreover, this implies 
$$
V_\lambda({\Bbb Z})=\bigoplus_{\mu\le \lambda} V_\lambda({\Bbb Z},\mu),
$$
where 
$$
V_\lambda({\Bbb Z},\mu)=\bigoplus_{{\bf n}\in{\Bbb N}_0^s\atop n_1\alpha_1+\cdots+n_s\alpha_s=\lambda-\mu} {\Bbb Z} X_-^{\bf n} {\rm v}_\lambda
=V_\lambda(\mu)\cap V_\lambda({\Bbb Z})\leqno(5)
$$ 
is the weight $\mu$ subspace of $V_\lambda({\Bbb Z})$ (the last "=" in equation (5) follows from equations (3) and (4)); in particular, $V_\lambda({\Bbb Z},\mu)$ is a ${\Bbb Z}$-lattice in $V_\lambda(\mu)$. 
The set of non-trivial weights in $V_\lambda$ is given as $\{\mu\in\Gamma_{\rm sc},\,\mu\le \lambda\}$ and the relative height function extends to a mapping
$$
{\rm ht}_\lambda:\,\{\mu\in\Gamma_{\rm sc},\,\mu\le\lambda\}\rightarrow{\Bbb N}_0
$$
(${\rm ht}_\lambda(\mu)={\rm ht}(\lambda-\mu)$; note that $\mu\le\lambda$ implies $\lambda-\mu\in\Gamma_{\rm ad}$). 
For any non-negative integer $r$ we define the following ${\Bbb Z}$-submodule of $V_\lambda$:
$$
V_\lambda({\Bbb Z},r)=\bigoplus_{\mu\le \lambda \atop 0\le {\rm ht}_\lambda(\mu)\le r} p^{r-{\rm ht}_\lambda(\mu)} V_\lambda({\Bbb Z},\mu)\oplus\bigoplus_{\mu\le \lambda \atop {\rm ht}_\lambda(\mu)> r} V_\lambda({\Bbb Z},\mu).
$$
Clearly, $V_\lambda({\Bbb Z},r)\le V_\lambda({\Bbb Z})$ is a ${\cal U}_{\Bbb Z}^-$-submodule.

\bigskip

Essentially the same proof as the one of (1.4) Lemma yields the following

\bigskip

{\bf Lemma. }{\it $V_\lambda({\Bbb Z},r)$ is an ${\cal S}$-invariant subspace of $V_\lambda({\Bbb Z})$.}



\bigskip

{\bf (2.2) Truncation of Verma modules. } We relate the truncation of irreducible, finite dimensional representations of ${\mathfrak g}$ to certain quotients 
("truncations") of Verma modules. This will be used in section 3. We assume that $\lambda$ is an integral and dominant weight, hence, the 
weight of any non-trivial weight space appearing in $V_\lambda$ is integral. For all simple roots $\alpha$ we set
$$
m_\alpha=\lambda(h_\alpha),
$$
hence, $\lambda=\sum_{\alpha\in \Delta} m_\alpha \omega_\alpha$. We denote by
$$
U_\lambda=\langle  x_{-\alpha}^{m_\alpha+1}{\rm v}_\lambda,\,\alpha\in\Delta \rangle_{\cal U}=\sum_{\alpha\in\Delta} {\cal U} x_{-\alpha}^{m_\alpha+1} {\rm v}_\lambda
$$
the ${\cal U}$-submodule of $V_\lambda$, which is generated by the elements $x_{-\alpha}^{m_\alpha+1}{\rm v}_\lambda$, $\alpha\in\Delta$. We claim that 
$$
U_\lambda=\langle  x_{-\alpha}^{m_\alpha+1}{\rm v}_\lambda,\,\alpha\in\Delta \rangle_{{\cal U}^-}=\sum_{\alpha\in\Delta} {\cal U}^- x_{-\alpha}^{m_\alpha+1} {\rm v}_\lambda.\leqno(6)
$$
In fact, the inclusion "$\supseteq$" is obvious. To prove the reverse inclusion, we let $\alpha\in \Delta$ and we denote by $S_\alpha\le {\mathfrak g}$ the Lie subalgebra isomorphic to ${\rm sl}_2$ which is attached to $\alpha$ as in [Hu], Proposition 8.3 (f), p. 37.
Applying the representation theory of ${\rm sl}_2$ to the ${\rm sl}_2({\Bbb C})$-module $S_\alpha{\rm v}_\lambda$ we find $x_\alpha x_{-\alpha}^{m_\alpha+1}{\rm v}_\lambda=0$.
Moreover $x_\beta x_{-\alpha}^{m_\alpha+1}{\rm v}_\lambda=0$ for any positive root $\beta\not=\alpha$ because $\beta-(m_\alpha+1)\alpha$ is not a sum of negative roots (note that $\alpha$ is simple). 
Since any $X\in{\cal U}$ is a linear combination of terms $X_-^{\bf a}H^{\bf b}X_+^{\bf c}$ we deduce that
$$
U_\lambda=\sum_{\alpha\in\Delta} {\cal U} x_{-\alpha}^{m_\alpha+1} {\rm v}_\lambda=\sum_{\alpha\in\Delta} {\cal U}^- x_{-\alpha}^{m_\alpha+1} {\rm v}_\lambda.
$$ 
which is the claim.

\medskip

Using that ${\cal U}^-=\bigoplus_{\bf n} {\Bbb C} X_-^{\bf n}$ and equation (6) we see that $U_\lambda$ is the sum of its weight spaces
$$
U_\lambda=\bigoplus_\mu U_\lambda(\mu),
$$
where $U_\lambda(\mu)=\sum_{\alpha,{\bf n}}{\Bbb C} X_-^{\bf n}x_{-\alpha}^{m_\alpha+1}{\rm v}_\lambda$ and $\alpha\in \Delta$, ${\bf n}=(n_i)\in{\Bbb N}_0^s$ run over all elements such that $\lambda-(m_\alpha+1)\alpha-\sum_i n_i\alpha_i=\mu$ (recall that $\alpha_1,\ldots,\alpha_s$ are the positive roots).
%
We set $U_\lambda({\Bbb Z})=U_\lambda\cap V_\lambda({\Bbb Z})$ and for any weight $\mu$ we put
$$
U_\lambda({\Bbb Z},\mu)=U_\lambda(\mu)\cap V_\lambda({\Bbb Z})=V_\lambda(\mu)\cap U_\lambda({\Bbb Z}) \quad(\le V_\lambda({\Bbb Z},\mu)).
$$
$U_\lambda({\Bbb Z})$ is a ${\cal U}_{\Bbb Z}$-module
and we obtain
$$
U_\lambda({\Bbb Z})=\bigoplus_\mu U_\lambda({\Bbb Z},\mu).
$$

\medskip

Theorem 21.4 in [Hu], p. 115 implies that there is a surjective morphism of ${\cal U}$-modules
$$
\varphi:\,V_\lambda\rightarrow L_\lambda,
$$
which sends $X {\rm v}_\lambda$ to $X v_\lambda$,  $X\in{\cal U}$, and has kernel $U_\lambda$. 
The restriction of $\varphi$ induces a surjective morphism of ${\cal U}_{\Bbb Z}$-modules
$$
\varphi_{\Bbb Z}:\,V_\lambda({\Bbb Z})\rightarrow L_\lambda({\Bbb Z})\leqno(7)
$$
which has kernel $U_\lambda({\Bbb Z})$. The morphism $\varphi_{\Bbb Z}$ induces a map
$$
\varphi_r:\,V_\lambda({\Bbb Z})\rightarrow L_\lambda({\Bbb Z})/L_\lambda({\Bbb Z},r)
$$
which sends ${\rm v}_\lambda$ to $\varphi_{\Bbb Z}({\rm v}_\lambda)+L_\lambda({\Bbb Z},r)=v_\lambda+L_\lambda({\Bbb Z},r)$. Since $L_\lambda({\Bbb Z},r)$ is an ${\cal S}$-module, $\varphi_r$ is a morphism of ${\cal S}$-modules.

\bigskip

{\bf Lemma. }{\it The kernel of $\varphi_r$ equals $U_\lambda({\Bbb Z})+V_\lambda({\Bbb Z},r)$. In particular, $U_\lambda({\Bbb Z})+V_\lambda({\Bbb Z},r)$ is an ${\cal S}$-module and $\varphi_r$ induces an isomorphism of ${\cal S}$-modules
$$
\frac{V_\lambda({\Bbb Z})}{U_\lambda({\Bbb Z})+V_\lambda({\Bbb Z},r)}\rightarrow \frac{L_\lambda({\Bbb Z})}{L_\lambda({\Bbb Z},r)}={\bf L}^{[r]}_\lambda({\Bbb Z}).
$$

}

\medskip

{\it Proof. } Since $\varphi$ is a morphism of ${\cal U}$-modules it respects weight spaces and since $\varphi$ maps $V_\lambda({\Bbb Z})$ to $L_\lambda({\Bbb Z})$ we see that $\varphi$ maps $V_\lambda({\Bbb Z},\nu)$ to $L_\lambda({\Bbb Z},\nu)$. We denote by $\varphi(\nu)$ the restriction of $\varphi$ to $V_\lambda({\Bbb Z},\nu)$. We claim that for any (integral) weight $\nu$ the sequence
$$
0\rightarrow U_\lambda({\Bbb Z},\nu)\stackrel{"\subseteq"}{\rightarrow} V_\lambda({\Bbb Z},\nu) \stackrel{\varphi(\nu)}{\rightarrow} L_\lambda({\Bbb Z},\nu)\rightarrow 0
$$
is exact. We first show surjectivity. Let $v\in L_\lambda({\Bbb Z},\nu)$. By the surjectivity of $\varphi_{\Bbb Z}$ there is $u\in V_\lambda({\Bbb Z})$ such 
that $\varphi(u)=v$. We write $u=\sum_\mu u_\mu$ with $u_\mu\in V_\lambda({\Bbb Z},\mu)$. Since $H \varphi(u)=\nu(H)\varphi(u)$ for all $H\in{\cal U}^o_{\Bbb Z}$ we obtain $Hu-\nu(H)u\in {\rm ker}\,\varphi_{\Bbb Z}=U_\lambda({\Bbb Z})$. 
Thus,
$$
\sum_\mu (\mu(H)-\nu(H)) u_\mu \in U_\lambda({\Bbb Z}).
$$
Since $U_\lambda({\Bbb Z})$ is the direct sum of its weight spaces $U_\lambda({\Bbb Z},\mu)$, this implies that
$(\mu(H)-\nu(H)) u_\mu\in U_\lambda({\Bbb Z},\mu)$ for all $\mu$. Since furthermore, $U_\lambda({\Bbb Z},\mu)\le V_\lambda({\Bbb Z})$ is saturated 
this implies $u_\mu\in U_\lambda({\Bbb Z},\mu)$ for all $\mu\not=\nu$. 
Hence, $\sum_{\mu\not=\nu} u_\mu\in U_\lambda({\Bbb Z})={\rm ker}\,\varphi_{\Bbb Z}$ and we obtain
$\varphi(u_\nu)=\varphi(u)=v$, which is the surjectivity of $\varphi(\nu)$.
The exactness of the sequence then is immediate because ${\rm ker}\,\varphi(\nu)= V_\lambda({\Bbb Z},\nu)\cap U_\lambda=U_\lambda({\Bbb Z},\nu)$.

\medskip

The surjectivity of $\varphi(\mu)$ implies that $\varphi(p^m \,V_\lambda({\Bbb Z},\mu))=p^m \,L_\lambda({\Bbb Z},\mu)$, hence,
$$
\varphi(V_\lambda({\Bbb Z},r))=L_\lambda({\Bbb Z},r).
$$ 
Since $U_\lambda({\Bbb Z})$ is the kernel of $\varphi_{\Bbb Z}$ this yields
$$
{\rm ker}\,\varphi_r=\varphi^{-1}_{\Bbb Z}(L_\lambda({\Bbb Z},r))=V_\lambda({\Bbb Z},r)+U_\lambda({\Bbb Z}).
$$
Hence, the proof is complete.

\bigskip

We set 
$$
\hat{V}_{\lambda}(r)=\bigoplus_{{\rm ht}_{\lambda}(\mu)>r} V_{\lambda}(\mu)\qquad {\rm and}\qquad\hat{L}_{\lambda}(r)=\bigoplus_{{\rm ht}_{\lambda}(\mu)>r} L_{\lambda}(\mu).
$$
We note that $\hat{V}_\lambda(r)$ and $\hat{L}_\lambda(r)$ are ${\cal U}^-$-modules.

\bigskip

{\bf Corollary. } {\it $(U_\lambda+\hat{V}_\lambda(r))\cap V_\lambda({\Bbb Z})\subseteq U_\lambda({\Bbb Z})+V_\lambda({\Bbb Z},r)$.

}

\medskip

{\it Proof. } Since the morphism of ${\cal U}$-modules $\varphi:\,V_{\lambda}\rightarrow L_{\lambda}$ has $U_{\lambda}$ in its kernel and respects weights spaces, it maps $U_\lambda+\hat{V}_{\lambda}(r)$ to $\hat{L}_{\lambda}(r)$. Since $\varphi(V_{\lambda}({\Bbb Z}))=L_{\lambda}({\Bbb Z})$ we obtain
$$
\varphi((U_{\lambda}+\hat{V}_{\lambda}(r))\cap V_{\lambda}({\Bbb Z}))\subseteq \hat{L}_{\lambda}(r)\cap L_{\lambda}({\Bbb Z})
=\bigoplus_{{\rm ht}_{\lambda}(\mu)>r} L_{\lambda}({\Bbb Z},\mu)
\subseteq L_{\lambda}({\Bbb Z},r).
$$
Hence, $(U_{\lambda}+\hat{V}_{\lambda}(r))\cap V_{\lambda}({\Bbb Z})$ is contained in the kernel of $\varphi_r$ which is identical with $U_{\lambda}({\Bbb Z})+V_{\lambda}({\Bbb Z},r)$ by the lemma. This completes the proof of the Corollary.

\bigskip

\section{Local constancy of the Truncated Representation}

{\bf (3.1) Some auxiliary Lemmas. }  We collect some auxiliary Lemmas which will be needed in the proofs in section (3.2) of the local constancy of the truncated Verma module. For any ${\bf n}=(n_i)_i\in{\Bbb N}_0^s$ and any ${\bf m}=(m_\gamma)_{\gamma\in\Delta}\in{\Bbb N}_0^\Delta$ we set
$$
X_-^{\cdot\bf n}=\prod_{i=1}^s x_{\alpha_i}^{n_i},\quad X_+^{\cdot\bf n}=\prod_{i=1}^s x_{-\alpha_i}^{n_i},\quad \mbox{and}\quad H^{\cdot\bf m}=\prod_{\gamma\in\Delta} h_\gamma^{m_\gamma}.
$$
Moreover, we define the length of ${\bf n}=(n_i)_i\in{\Bbb N}_0^s$ as $\ell({\bf n})=\sum_i n_i$.

\medskip

{\bf (3.1.1) Lemma. } {\it Let $\alpha\in\Phi^+$ be any positive root and let ${\bf n}=(n_i)_i\in{\Bbb N}_0^s$.  Then
$$
x_\alpha X_-^{\cdot \bf n}=\sum_{{\bf a},{\bf b}\in{\Bbb N}_0^s} \zeta_{{\bf a},{\bf b}} X_-^{\cdot \bf a}X_+^{\cdot \bf b}
+\sum_{{\bf a}\in{\Bbb N}_0^s,\gamma\in\Phi^+} \zeta_{{\bf a},\gamma} X_-^{\cdot \bf a} h_\gamma,\leqno(1)
$$
with certain $\zeta_{{\bf a},{\bf b}},\zeta_{{\bf a},\gamma}\in{\Bbb C}$.

}

{\it Proof. } We denote by $i$ the smallest index such that $n_i\not=0$ and we set $\beta=\alpha_i$ ($\alpha_1,\ldots,\alpha_s$ is the ordering of $\Phi^+$;  cf. section (1.1)). Hence, we may write 
$X_-^{\cdot\bf n}=x_{-\beta} X_-^{\cdot{\bf n}'}$, where $\ell({\bf n}')=\ell({\bf n})-1$ and we obtain

\begin{eqnarray*}
(2)\qquad x_\alpha X_-^{\cdot\bf n}&=&x_\alpha x_{-{\beta}} X_-^{\cdot{\bf n}'}\\
&=&x_{-{\beta}} x_\alpha X_-^{\cdot{\bf n}'}+[x_\alpha,x_{-{\beta}}] X_-^{\cdot{\bf n}'}.\\
\end{eqnarray*}

We want to show that $[x_\alpha,x_{-{\beta}}] X_-^{\cdot{\bf n}'}$ is of the form as on the right hand side of equation (1), i.e. 
$$
[x_\alpha,x_{-{\beta}}] X_-^{\cdot{\bf n}'}=\sum_{{\bf a},{\bf b}\in{\Bbb N}_0^s} \zeta'_{{\bf a},{\bf b}} X_-^{\cdot \bf a}X_+^{\cdot \bf b}
+\sum_{{\bf a}\in{\Bbb N}_0^s,\gamma\in\Phi^+} \zeta'_{{\bf a},\gamma} X_-^{\cdot \bf a} h_\gamma\leqno(3) 
$$
for certain $\zeta'_{{\bf a},{\bf b}},\zeta_{{\bf a},\gamma}\in{\Bbb C}$. To this end we note that there are four possibilities for the difference $\alpha-\beta$:
\begin{itemize}
\item Case 1. $\alpha-\beta$ is a positive root. In this case Theorem 25.2 in [Hu], p. 147 implies that $[x_\alpha,x_{-\beta}]=\zeta x_{\alpha-\beta}$ for some $\zeta\in{\Bbb Z}$  

\item Case 2. $\alpha-\beta\not=0$ is no root, hence, $[x_\alpha,x_{-\beta}]=0$

\item Case 3. $\alpha-\beta=0$. In this case Theorem 25.2 in [Hu], p. 147 implies that $[x_\alpha,x_{-\beta}]=\sum_{\gamma\in\Delta} \zeta_\gamma h_\gamma$ for certain $\zeta_\gamma\in{\Bbb Z}$ 

\item Case 4. $\alpha-\beta$ is a negative root. Again, in this case Theorem 25.2 in [Hu], p. 147 implies that $[x_\alpha,x_{-\beta}]=\zeta x_{\alpha-\beta}$ for some $\zeta\in{\Bbb Z}$. 
\end{itemize}

In the first case using an induction over the length of ${\bf n}'$, equation (1) (with the positive root $\alpha$ replaced by the positive root 
$\alpha-\beta$ and ${\bf n}$ replaced by ${\bf n}'$) shows that  $[x_\alpha,x_{-\beta}] X_-^{\cdot{\bf n}'}=\zeta x_{\alpha-\beta} X_-^{\cdot{\bf n}'}$ is of the 
form (3). In the second case $[x_\alpha,x_{-\beta}] X_-^{\cdot{\bf n}'}$ vanishes, hence equation (3) holds trivially. In the third case for any positive root 
$\gamma$ we have $h_\alpha x_{-\gamma}=x_{-\gamma}h_\alpha+[h_\alpha,x_{-\gamma}]=x_{-\gamma}h_\alpha-\gamma(h_\alpha) x_{-\gamma}$, where $\gamma(h_\alpha)\in{\Bbb Z}$. 
A simple induction over the length of ${\bf n}'$ therefore shows that 
$$
h_\alpha X_-^{\cdot{\bf n'}}=X_-^{\cdot{\bf n}'}h_\alpha+\zeta X_-^{\cdot{\bf n}'}
$$
for some $\zeta\in{\Bbb C}$, hence, equation (3) holds. In the fourth case Lemma 5.14 in [Ha], p. 135 implies that  $[x_\alpha,x_{-\beta}] X_-^{\cdot{\bf n}'}=\zeta x_{\alpha-\beta}X_-^{\cdot{\bf n}'}$ can be written
$$
[x_\alpha,x_{-\beta}] X_-^{\cdot{\bf n}'}=\sum_{{\bf a}\in{\Bbb N}_0^s} \zeta_{\bf a} X_-^{\cdot\bf a}
$$
for some $\zeta_{\bf a}\in{\Bbb C}$. Thus, equation (3) is proven and we have shown that
$$
x_\alpha X_-^{\cdot\bf n}=x_{-\beta} x_{\alpha} X_-^{\cdot{\bf n}'}+\mbox{terms as in the right hand side of eq. (1)},
$$ 
where now $\ell({\bf n}')< \ell({\bf n})$. Repeating the above computation we find
$$
x_\alpha X_-^{\cdot {\bf n}'}=x_{-\beta'} x_{\alpha} X_-^{\cdot{\bf n}''}+\mbox{terms as in the right hand side of eq. (1)},
$$ 
where now $\ell({\bf n}'')< \ell({\bf n}')$ and $\beta'=\alpha_i$, $i$ the smallest index such that $n'_i\not=0$. Proceeding in this way we finally obtain that
$$
x_\alpha X_-^{\cdot \bf n}=X_-^{\cdot {\bf n}} x_{-\alpha} + \mbox{terms as in the right hand side of eq. (1)}.
$$
This is the claim and the proof of the Lemma is complete.

\medskip

We define the length of ${\bf m}=(m_\gamma)_{\gamma}\in{\Bbb N}_0^\Delta$ as $\ell({\bf m})=\sum_{\gamma\in\Delta} m_\gamma$.

\medskip

{\bf (3.1.2) Lemma. }{\it For any positive root $\alpha$ and any $k\in{\Bbb N}$ 
$$
x_\alpha^k X_-^{\cdot{\bf n}}=\sum_{{\bf a},{\bf c}\in{\Bbb N}_0^s,\,{\bf b}\in{\Bbb N}_0^\Delta\atop \ell({\bf b})\le k} \zeta_{{\bf a},{\bf b},{\bf c}} X_-^{\cdot\bf a} H^{\cdot{\bf b}} X_+^{\cdot\bf c}.
$$

 }

{\it Proof. } We use induction on $k$. The case $k=1$ is immediate by the preceeding Lemma. We write using the induction hypothesis
$$
x_\alpha^{k+1}X_-^{\cdot{\bf n}}=x_\alpha x_\alpha^k X_-^{\cdot{\bf n}}=x_\alpha\big(  \sum_{{\bf a},{\bf c}\in{\Bbb N}_0^s,\,{\bf b}\in{\Bbb N}_0^\Delta\atop \ell({\bf b})\le k} \zeta_{{\bf a},{\bf b},{\bf c}} X_-^{\cdot \bf a} H^{\cdot\bf b} X_+^{\cdot \bf c}   \big).
$$
Using the expression for $x_\alpha X_-^{\cdot\bf a}$ from the preceding Lemma we obtain
$$
x_\alpha^{k+1}X_-^{\cdot{\bf n}}
=\sum_{{\bf a},{\bf c}\in{\Bbb N}_0^s,\,{\bf b}\in{\Bbb N}_0^\Delta\atop \ell({\bf b})\le k} \zeta_{{\bf a},{\bf b},{\bf c}} \sum_{{\bf a}',{\bf b}'} \zeta_{{\bf a}',{\bf b}'} 
X_-^{\cdot{\bf a}'} X_+^{\cdot{\bf b}'} H^{\cdot\bf b} X_+^{\cdot \bf c} +\sum_{{\bf a}',\gamma\in\Phi^+} X_-^{\cdot{\bf a}'} h_\gamma H^{\cdot\bf b} X_+^{\cdot \bf c}.\leqno(4)
$$
Taking into account that for any positive root $\gamma$
$$
X_+^{\cdot{\bf b}'}h_\gamma=h_\gamma X_+^{\cdot{\bf b}'}+\zeta X_+^{\cdot{\bf b}'}
$$
for some $\zeta\in{\Bbb C}$ we see that $X_+^{\cdot {\bf b}'} H^{\cdot \bf b}$ can be written as a sum $\sum_{{\bf n},\,\ell({\bf n})\le \ell({\bf b})} \zeta_{\bf n} H^{\cdot \bf n} X_+^{\cdot {\bf b}'}$; hence, the summation over ${\bf a}',{\bf b}'$ in equation (4) is of the form as claimed in the Lemma. Since ${\cal U}^o$ is 
abelian we see that the summation over ${\bf a}',\gamma$ in equation (4) also is of the form as claimed in the lemma. This completes the proof of the 
Lemma.

\medskip

{\bf (3.1.3) Lemma. }{\it For any positive root $\alpha\in\Phi^+$ and any $k\in{\Bbb N}$ 
$$
\frac{x_\alpha^k}{k!} X_-^{{\bf n}}=\sum_{{\bf a},{\bf b},{\bf c}\in{\Bbb N}_0^s\atop \ell({\bf b})\le k} \zeta_{{\bf a},{\bf b},{\bf c}} X_-^{\bf a} H^{{\bf b}} X_+^{\bf c}.
$$
with certain $\zeta_{{\bf a},{\bf b},{\bf c}}\in{\Bbb Z}$.
}

{\it Proof. } The ${\Bbb C}$-linear span $\langle H^{\cdot{\bf n}},\,\ell({\bf n})\le k \rangle_{\Bbb C}$ contains the ${\Bbb C}$-linear span $\langle H^{{\bf n}},\,\ell({\bf n})\le k \rangle_{\Bbb C}$ and since $H^{\bf n}$ has "leading monomial" $\prod_{\gamma\in\Delta} h_\gamma^{n_\gamma}=H^{\cdot\bf n}$ we find that
$$
\langle H^{\cdot{\bf n}},\,\ell({\bf n})\le k \rangle_{\Bbb C}=\langle H^{{\bf n}},\,\ell({\bf n})\le k \rangle_{\Bbb C}.
$$
Thus, using the preceding Lemma we deduce that we can write
$$
\frac{x_\alpha^k}{k!} X_-^{{\bf n}}=\sum_{{\bf a},{\bf b},{\bf c}\in{\Bbb N}_0^s\atop \ell({\bf b})\le k} \zeta_{{\bf a},{\bf b},{\bf c}} X_-^{\bf a} H^{{\bf b}} X_+^{\bf c}
$$
with certain $\zeta_{{\bf a},{\bf b},{\bf c}}\in{\Bbb C}$. On the other hand, $\frac{x_\alpha^k}{k!} X_-^{{\bf n}}$ is contained in the ${\Bbb Z}$-lattice ${\cal U}_{\Bbb Z}$ of ${\cal U}$ and since
${\cal U}_{\Bbb Z}$ is a free ${\Bbb Z}$-module with basis $X_-^{\bf a}H^{\bf b}X_+^{\bf c}$ we deduce that $\zeta_{{\bf a},{\bf b},{\bf c}}\in{\Bbb Z}$. Thus, the Lemma is proven.

\medskip

\bigskip

{\bf (3.2) Local Constancy of truncated Verma modules. } The following Proposition shows that the truncations of two Verma modules $V_\lambda({\Bbb Z})$ and $V_{\lambda'}({\Bbb Z})$ are isomorphic if the highest weights $\lambda$ and $\lambda'$ are sufficiently close in the $p$-adic sense.

\medskip

{\bf Proposition. }{\it Let $r\in {\Bbb N}$, let ${\cal T}\subseteq\Delta$ be a subset and let $\lambda,\lambda'\in\Gamma_{\rm sc}$ be integral and dominant weights. We set $m_\alpha=\lambda(h_\alpha)$, $m_\alpha'=\lambda'(h_\alpha)$, $\alpha\in \Delta$. If

\begin{itemize}

\item $m_\alpha=m_\alpha'$ for all $\alpha\in \Delta-{\cal T}$

\item $m_\alpha>r$, $m_\alpha'>r$ for all $\alpha\in {\cal T}$

\item $\lambda\equiv \lambda'\pmod{p^{\lceil\frac{p}{p-1}r\rceil}\Gamma_{\rm sc}}$ 

\end{itemize}

then 
$$
\frac{V_\lambda({\Bbb Z})}{U_\lambda({\Bbb Z})+V_\lambda({\Bbb Z},r)}\cong  \frac{V_{\lambda'}({\Bbb Z})}{U_{\lambda'}({\Bbb Z})+V_{\lambda'}({\Bbb Z},r)}
$$
as ${\cal S}$-modules.

}

\medskip

{\it Proof. } As before, we write $V_\lambda({\Bbb Z})={\cal U}^-_{\Bbb Z}{\rm v}_\lambda$ and $V_{\lambda'}({\Bbb Z})={\cal U}_{\Bbb Z}^-
{\rm v}_{\lambda'}$. Since $V_\lambda({\Bbb Z})$ and $V_{\lambda'}({\Bbb Z})$ are free ${\cal U}_{\Bbb Z}^-$-modules, there is a uniquely determined isomorphism of ${\cal U}_{\Bbb Z}^-$-modules
$$
\Phi:\,V_\lambda({\Bbb Z})\rightarrow V_{\lambda'}({\Bbb Z}),\leqno(5)
$$
which sends ${\rm v}_\lambda$ to ${\rm v}_{\lambda'}$, hence, $\Phi(X{\rm v}_\lambda)=X{\rm v}_{\lambda'}$ for all 
$X\in {\cal U}_{\Bbb Z}^-$. We claim that
$$
\Phi(U_\lambda({\Bbb Z})+V_\lambda({\Bbb Z},r))= U_{\lambda'}({\Bbb Z})+V_{\lambda'}({\Bbb Z},r).\leqno(6)
$$
We will prove the inclusion "$\subseteq$"; the proof of the reverse inclusion follows in the same way by considering the inverse isomorphism $\Phi^{-1}:\,V_{\lambda'}({\Bbb Z})\rightarrow V_\lambda({\Bbb Z})$ which sends ${\rm v}_{\lambda'}$ to ${\rm v}_\lambda$. We will seperately show that
$$
\Phi(V_\lambda({\Bbb Z},r))\subseteq V_{\lambda'}({\Bbb Z},r) \leqno(6a)
$$
and
$$
\Phi(U_\lambda({\Bbb Z}))\subseteq  U_{\lambda'}({\Bbb Z})+V_{\lambda'}({\Bbb Z},r). \leqno(6b)
$$ 
We begin with the proof of equation (6a). We select an arbitrary (integral) weight
$\mu\le \lambda$ and an arbitrary vector $v\in V_\lambda({\Bbb Z},\mu)$. In view of equation (5) in section (2.1), $v$ can be written $v=\sum_{{\bf n}} c_{\bf
n}X_-^{\bf n} {\rm v}_\lambda$, where $c_{\bf n}\in{\Bbb Z}$ and the sum runs over all ${\bf n}=(n_i)\in {\Bbb N}_0^s$ satisfying 
$\sum_{i=1}^s n_i \alpha_i=\lambda-\mu$. Applying $\Phi$ yields $\Phi(v)=\sum_{\bf n} c_{\bf n} X_-^{\bf n} {\rm v}_{\lambda'}$, hence,
$\Phi(v)$ is contained in $V_{\lambda'}({\Bbb Z},\lambda'-(\lambda-\mu))$ and we have shown that $\Phi(V_\lambda({\Bbb Z},\mu))\subseteq V_{\lambda'}({\Bbb Z},\lambda'-(\lambda-\mu))$. Since
$$
{\rm ht}_{\lambda'}(\lambda'-(\lambda-\mu))=\lambda-\mu={\rm ht}_\lambda(\lambda-(\lambda-\mu))={\rm ht}_\lambda(\mu)
$$
we thus have shown that
$$
\Phi(V_\lambda({\Bbb Z},\mu))\subseteq V_{\lambda'}({\Bbb Z},\mu')\leqno(7)
$$
for some integral weight $\mu'\le\lambda'$ satisfying ${\rm ht}_{\lambda'}(\mu')={\rm ht}_\lambda(\mu)$ ($\mu'=\lambda'-(\lambda-\mu)$). 
By definition of $V_\lambda({\Bbb Z},r)$ this implies that equation (6a) holds. To prove equation (6b) we recall that $U_\lambda$ is generated as ${\cal U}^-$-module by $x_{-\alpha}^{m_\alpha+1}{\rm v}_\lambda$, $\alpha\in \Delta$.  For all $\alpha\in\Delta-{\cal T}$ we find 
$$
\Phi(x_{-\alpha}^{m_\alpha+1}{\rm v}_\lambda)=x_{-\alpha}^{m_\alpha+1}{\rm v}_{\lambda'}=x_{-\alpha}^{m_\alpha'+1}{\rm v}_{\lambda'}\in U_{\lambda'}.
$$

\medskip

We let $\alpha\in{\cal T}$. Since $\Phi(x_{-\alpha}^{m_\alpha+1}{\rm v}_\lambda)=x_{-\alpha}^{m_\alpha+1}{\rm v}_{\lambda'}$ has weight $\lambda'-(m_\alpha+1)\alpha$ and $m_\alpha> r$ for $\alpha\in {\cal T}$ we find
$$
\Phi(x_{-\alpha}^{m_\alpha+1}{\rm v}_\lambda)\in \hat{V}_{\lambda'}(r).
$$
Altogether, we obtain for all $\alpha\in\Delta$ that $\Phi(x_{-\alpha}^{m_\alpha+1}{\rm v}_\lambda)\in U_{\lambda'}+\hat{V}_{\lambda'}(r)$, hence, $\Phi(U_\lambda)\subseteq U_{\lambda'}+\hat{V}_{\lambda'}(r)$ (note that $U_{\lambda'}+\hat{V}_{\lambda'}(r)$ is a ${\cal U}^-$-module).
Since $\Phi(V_\lambda({\Bbb Z}))=V_{\lambda'}({\Bbb Z})$  by definition of $\Phi$, we obtain
$$
\Phi(U_\lambda({\Bbb Z}))\subseteq (U_{\lambda'}+\hat{V}_{\lambda'}(r))\cap V_{\lambda'}({\Bbb Z}).
$$
The right hand side of the above equation is contained in $U_{\lambda'}({\Bbb Z})+V_{\lambda'}({\Bbb Z},r)$ by (2.2) Corollary, hence, equation (6b) follows.

\medskip

Equation (6) implies that $\Phi$ induces an isomorphism of ${\cal U}_{\Bbb Z}^-$-modules
$$
\Phi:\,V_\lambda({\Bbb Z})/(U_\lambda({\Bbb Z})+V_\lambda({\Bbb Z},r))\rightarrow V_ {\lambda'}({\Bbb Z})/(U_{\lambda'}({\Bbb Z})+V_{\lambda'}({\Bbb Z},r)),
$$
i.e. $\Phi$ commutes with all $\frac{x_{-\alpha}^n}{n!}$, $\alpha\in\Phi^+$. 
It remains to show that $\Phi$ commutes with the action of all generators of type "$p^{n {\rm ht}(\alpha)}\frac{x_\alpha^n}{n!}$" of ${\cal S}$, i.e. 
$$
\Phi(t^n\frac{x_\alpha^n}{n!} v)=t^n\frac{x_\alpha^n}{n!} \Phi(v)\leqno(8)
$$ 
for all $\alpha\in \Phi^+$, $n\in{\Bbb N}_0$, $t\in p^{{\rm ht}(\alpha)}{\Bbb Z}$ and $v\in V_\lambda({\Bbb Z})/(U_\lambda({\Bbb Z})+V_\lambda({\Bbb Z},r))$
(this implies that $\Phi$ is an isomorphism of ${\cal S}$-modules as claimed; note that 
$V_\lambda({\Bbb Z})/(U_\lambda({\Bbb Z})+V_\lambda({\Bbb Z},r))$ and $V_ {\lambda'}({\Bbb Z})/(U_{\lambda'}({\Bbb Z})+V_{\lambda'}({\Bbb Z},r))$ are ${\cal S}$-modules; cf. (2.2) Lemma). 
Since $\Phi$ is ${\Bbb Z}$-linear we may assume that $v$ is of the form 
$$
v=X^{\bf n}_- {\rm v}_\lambda+(U_\lambda({\Bbb Z})+V_\lambda({\Bbb Z},r))=X_-^{\bf n}\bar{\rm v}_\lambda\qquad(\bar{\rm v}_\lambda={\rm v}_\lambda+(U_\lambda({\Bbb Z})+V_\lambda({\Bbb Z},r)))
$$ 
for some ${\bf n}\in {\Bbb N}_0^s$. 
%
%
%
%
Since $p$ divides $p^{{\rm ht}(\alpha)}$, it divides $t$ and since $p^r$ annihilates $V_\lambda({\Bbb Z})/V_\lambda({\Bbb Z},r)$ and $V_{\lambda'}({\Bbb Z})/V_{\lambda'}({\Bbb Z},r)$, we see that both sides of equation (8) vanish if $n\ge r$. We therefore may assume that 
$$
n<r.
$$
Using (3.1.3) Lemma we can write
$$
\frac{x_\alpha^n}{n!} X_-^{\bf n}=\sum_{{\bf a},{\bf b},{\bf c} \atop \ell({\bf b})\le n } \zeta_{{\bf a},{\bf b},{\bf c}} 
X_-^{{\bf a}}H^{{\bf b}}X_+^{{\bf c}},
$$
where $\zeta_{{\bf a},{\bf b},{\bf c}}\in{\Bbb Z}$ and $\ell({\bf b})\le n<r$. 
%
%
%
%
%
%
Since $H^{\bf b}{\rm v}_\lambda=\lambda(H^{\bf b}){\rm v}_\lambda$ and since $X_+^{\bf n}$ annihilates ${\rm v}_\lambda$ we obtain 
$$
t^n\,\frac{x_\alpha^n}{n!} X_-^{\bf n}{\rm v}_\lambda=t^n\,\sum_{{\bf a},{\bf b} \atop \ell({\bf b})\le n } \zeta_{{\bf a},{\bf b},0} \lambda(H^{\bf b}) 
X_-^{{\bf a}}{\rm v}_\lambda.\leqno(9)
$$
Reducing equation (9) modulo $U_\lambda({\Bbb Z})+V_\lambda({\Bbb Z},r)$ and applying $\Phi$ we obtain
\begin{eqnarray*}
(10)\qquad\Phi(t^n\frac{x_\alpha^n}{n!}X_-^{\bf n} \bar{\rm v}_\lambda)&=&\Phi(t^n\,\sum_{{\bf a},{\bf b}} \zeta_{{\bf a},{\bf b},0} \,
\lambda(H^{\bf b}) X_-^{{\bf a}} \bar{\rm v}_\lambda)\\
&=&t^n\,\sum_{{\bf a},{\bf b}\atop \ell({\bf b})\le n} \zeta_{{\bf a},{\bf b},0} \,\lambda(H^{\bf b}) X_-^{{\bf a}} \bar{\rm v}_{\lambda'}
\end{eqnarray*}
(note that $t^n\,\frac{x_\alpha^n}{n!} X_-^{\bf n}$ and $\zeta_{{\bf a},{\bf b},0} \,\lambda(H^{\bf b}) X_-^{{\bf a}}$ are contained in ${\cal S}$, 
hence, these act well defined on $\bar{\rm v}_\lambda$, $\bar{\rm v}_{\lambda'}$). 
%
In the same way we find
$$
t^n\,\frac{x_\alpha^n}{n!} \Phi(X_-^{\bf n} \bar{\rm v}_\lambda)=t^n\,\frac{x_\alpha^n}{n!}  X_-^{\bf n} \bar{\rm v}_{\lambda'}
=t^n\,\sum_{{\bf a},{\bf b}\atop \ell({\bf b})\le n} \zeta_{{\bf a},{\bf b},0} \lambda'(H^{\bf b}) X_-^{{\bf a}}\bar{\rm v}_{\lambda'}.\leqno(11)
$$
For any weight $\mu\in\Gamma_{\rm sc}$ we know that $\mu(h_\beta)=\langle \mu,\beta\rangle\in{\Bbb Z}$ for all simple roots $\beta\in\Delta$. 
The congruence $\lambda\equiv\lambda'\pmod{p^{\lceil\frac{p}{p-1}r\rceil}\Gamma_{\rm sc}}$ therefore implies that
$$
\lambda(h_\beta)-\lambda'(h_\beta)\in p^{\lceil\frac{p}{p-1}r\rceil}{\Bbb Z}.\leqno(12)
$$
On the other hand, since $\ell({\bf b})<r$ we know that $b_\beta<r$ for all $\beta\in\Delta$. Hence, 
$v_p(b_{\beta}!)<\frac{r}{p-1}$ and together with equation (12) we find
$$
v_p\left({\lambda(h_\beta)\choose b_{\beta}}-{\lambda'(h_\beta)\choose b_{\beta}} \right)\ge r
$$
for all $\beta\in\Delta$, hence,
$$
\lambda(H^{\bf b}) \equiv \lambda'(H^{\bf b})\pmod{p^r}.
$$
Together with equations (10) and (11) and taking into account that $p^r$ annihilates $V_{\lambda'}({\Bbb Z})/V_{\lambda'}({\Bbb Z},r)$ this finally yields
$$
\Phi(t^n\frac{x_\alpha^n}{n!}X_-^{\bf n}\bar{\rm v}_\lambda)=t^n\,\frac{x_\alpha^n}{n!} \Phi(X_-^{\bf n} \bar{\rm v}_\lambda).
$$
Thus, $\Phi$ is ${\cal S}$-invariant and the proof of the Proposition is complete.

\bigskip



\bigskip

{\bf (3.3) Constancy of truncated finite dimensional irreducible representations. } (2.2) Lemma and (3.2) Proposition imply

\bigskip

{\bf Proposition. }{\it 
%
Let $r\in {\Bbb N}$, let ${\cal T}\subseteq\Delta$ be a subset and let $\lambda,\lambda'\in\Gamma_{\rm sc}$ be integral and dominant weights. We set $m_\alpha=\lambda(h_\alpha)$, $m_\alpha'=\lambda'(h_\alpha)$, $\alpha\in \Delta$. If

\begin{itemize}

\item $m_\alpha=m_\alpha'$ for all $\alpha\in \Delta-{\cal T}$

\item $m_\alpha>r$, $m_\alpha'>r$ for all $\alpha\in {\cal T}$

\item $\lambda\equiv \lambda'\pmod{p^{\lceil\frac{p}{p-1}r\rceil}\Gamma_{\rm sc}}$ 

\end{itemize}

then 
$$
{\bf L}^{[r]}_\lambda({\Bbb Z})\cong{\bf L}^{[r]}_{\lambda'}({\Bbb Z})
$$
as ${\cal S}$-modules.

}

\bigskip

\section{Representations of split semi simple groups}

{\bf (4.1) Chevalley groups. } We fix a dominant and integral weight $\lambda_0\in \Gamma_{\rm sc}$ such that the representation $(\rho_{\lambda_0},L_{\lambda_0})$ of ${\mathfrak g}$ is faithful (we can always assume this by omitting the simple factors in the kernel of $\rho_{\lambda_0}$). Let $R$ denote any 
${\Bbb Z}$-algebra. For any root $\alpha\in\Phi$ and any $t\in R$ we set
$$
x_\alpha(t)=\exp \rho_{\lambda_0}(t x_\alpha)\in {\rm SL}(L_{\lambda_0}(R)).
$$
We then denote by 
$$
{G}_{{\lambda_0},R}=\langle x_\alpha(t_\alpha),\;\alpha\in\Phi, \,t_\alpha\in R\, \rangle\:\le {\rm SL}(L_{\lambda_0}(R))\leqno(1)
$$ 
the {\it Chevalley group} attached to $\rho_{\lambda_0}$ and $R$, i.e. $G_{\lambda_0,R}$ is the subgroup of ${\rm
SL}(L_{\lambda_0}(R))$, which is generated by the $x_\alpha(t)$, $\alpha\in\Phi$, $t\in R$. Equation (1) even defines an algebraic group in the following 
sense: we choose a ${\Bbb Z}$-basis of $L_{\lambda_0}({\Bbb Z})$, which yields an identification
$$
{\rm SL}(L_{\lambda_0}(R))={\rm SL}_m(R)\qquad(m={\rm dim}\,L_{\lambda_0}).
$$
Hence, $G_{\lambda_0,R}\le {\rm SL}_m(R)$ and there is an algebraic group ${\bf G}_{\lambda_0}\le {{\bf SL}_m}_{/{\Bbb Q}}$, 
which is defined over ${\Bbb Q}$, such that 
$$
{\bf G}_{\lambda_0}(\bar{\Bbb Q})=G_{\lambda_0,\bar{\Bbb Q}}.
$$

\medskip

We denote by ${\bf SL}(L_\lambda({\Bbb Z}))$ the algebraic ${\Bbb Z}$-group, which is defined by ${\bf SL}(L_\lambda({\Bbb Z}))(R)={\bf SL}(L_\lambda(R))$, $R$ 
any ${\Bbb Z}$-algebra. The choice of a ${\Bbb Z}$-basis of $L_\lambda({\Bbb Z})$ yields an identification ${\bf SL}(L_\lambda({\Bbb Z}))={\bf SL}_m/{\Bbb Z}$, 
where $m={\rm dim}\,L_\lambda({\Bbb Z})$. The group ${\bf G}_{\lambda_0}$ then has a natural ${\Bbb Z}$-structure ${{\bf G}_{\lambda_0}}_{/{\Bbb Z}}$ (i.e. a ${\Bbb Z}$-form 
of the coordinate algebra ${\Bbb Q}[{\bf G}_{\lambda_0}]$), such that ${{\bf G}_{\lambda_0}}_{/{\Bbb Z}}$ embeds as a closed subscheme of 
${\bf SL}(L_\lambda({\Bbb Z}))={{\bf SL}_m}_{/{\Bbb Z}}$ (cf. [B], sec. 3.4, p. 18). Hence, if $R$ is any ${\Bbb Z}$-algebra, which is embedded in $\bar{\Bbb Q}$, 
then
$$
{\bf G}_{\lambda_0}(R)={\bf G}_{\lambda_0}(\bar{\Bbb Q})\cap {\rm SL}_m(R).
$$
In particular, we obtain
$$
G_{\lambda_0,R}\subseteq {\rm SL}(L_{\lambda_0}(R))\cap G_{\lambda_0,\bar{\Bbb Q}}
={\rm SL}_m(R)\cap {\bf G}_{\lambda_0}(\bar{\Bbb Q})
={\bf G}_{\lambda_0}(R).\leqno(2)
$$

\medskip

The group ${\bf G}_{\lambda_0}$ even is split over ${\Bbb Q}$ and it is well known that any semi-simple, ${\Bbb Q}$-split algebraic group 
${\bf G}$ is isomorphic to ${\bf G}_{\lambda_0}$ for some faithful representation $\rho_{\lambda_0}$ of the Lie algebra ${\mathfrak g}$ of ${\bf G}$.

\bigskip

{\bf (4.2) } For any prime $p\in{\Bbb N}$ we define the level subgroup
$$
K_*(p)\le{\bf G}_{\lambda_0,{\Bbb Z}_p}
$$ 
as the subgroup, which is generated by the elements $x_\alpha(t_\alpha)$, where $\alpha\in\Phi$ and 
$$
t_\alpha\in\left\{ \begin{array}{ccc} p^{{\rm ht}(\alpha)}{\Bbb Z}_p&\mbox{if}&\alpha\in\Phi^+\\ 
{\Bbb Z}_p&\mbox{if}&\alpha\in\Phi^-.
\end{array}\right.
$$
Equation (2) implies that 
$$
K_*(p)\le{\bf G}_{\lambda_0}({\Bbb Z}_p).
$$ 


\medskip

{\it Example. } We set ${\mathfrak g}={\mathfrak sl}_n$ and we denote by ${\mathfrak h}$ the Cartan subalgebra consisting of 
diagonal matrices in ${\mathfrak g}$. The roots of ${\mathfrak g}$ with respect to ${\mathfrak h}$ are $\alpha_{ij}$, $i\not=j$, 
where $\alpha_{ij}({\rm diag}(h_1,\ldots,h_n))=h_i-h_j$. The set $\Delta=\{\alpha_{i+1,i},\;i=1,\ldots,n-1\}$ is a basis for 
$\Phi$ and the root space ${\mathfrak g}(\alpha_{ij})$ is generated by the elementary matrix
$x_{\alpha_{ij}}=e_{ij}=(\delta_{ij})_{ij}$.
The height of a root is explicitly given as ${\rm ht}(\alpha_{ij})=i-j$. We fix the standard matrix representation
$\rho_0:\,{\mathfrak g}\rightarrow {\mathfrak sl}_n$ of ${\mathfrak g}$, i.e. $\rho_0(x)=x$ for all $x\in{\mathfrak g}$. An easy computation shows that
$$
e_{ij}(t):=x_{\alpha_{ij}}(t)=\exp \rho_0(te_{ij})=1+te_{ij}\leqno(3)
$$
for all $i\not=j$. This implies in particular that ${\bf G}_{\rho_0}={\bf SL}_n$. We make the following
$$
{\it Claim. }\quad K_*(p)=\{(\gamma_{ij})\in {\rm SL}_n({\Bbb Z}_p):\, \gamma_{ii}\in 1+{p}{\Bbb Z}_p\;\mbox{for all}\; i,\;\;\gamma_{ij}\in p^{i-j}{\Bbb Z}_p
\,\mbox{for all}\; i>j\}.\leqno(4)
$$

{\it Proof. } It is easily verified that the right hand side defines a subgroup of ${\rm SL}_n({\Bbb Z}_p)$ and equation (3) implies that 
all generators $x_{\alpha_{ij}}(t)=e_{ij}(t)$ of $K_*(p)$ (i.e. $t\in{\Bbb Z}_p$ if $i<j$ and $t\in p^{{\rm ht}(\alpha_{ij})}{\Bbb Z}_p=p^{i-j}{\Bbb Z}_p$ if $i>j$) 
are contained in the right hand side of equation (4). Hence, the inclusion "$\subseteq$" holds. To prove the reverse inclusion "$\supseteq$" let $\gamma$ be contained in the right hand side of equation (4). Using left multiplication by elements $x_{{\alpha}_{ij}}(t)$ 
with $i>j$ and $t\in p^{i-j}{\Bbb Z}_p$, i.e. by using row operations, we can transform $\gamma$ into an upper triangular matrix 
$\gamma'=(\gamma'_{ij})$ (note that $\gamma_{ij}\in p^{i-j}{\Bbb Z}_p$ if $i>j$). Since $\gamma'$ still has determinant equal to $1$ we know that $\prod_i \gamma_{ii}'=1$. The ${\rm SL}_2$-relation
$$
\Mat{t^{-1}}{}{}{t}=\Mat{1}{}{(t-1)t}{1}\Mat{1}{t^{-1}}{}{1}\Mat{1}{}{-(t-1)}{1}\Mat{1}{-1}{}{1}
$$
implies that $K_*(p)$ contains the matrices
$$
h_{{i,i+1}}(t)=\left(\begin{array}{ccccc}
\ddots&&&\\
&t&&\\
&&t^{-1}&\\
&&&\ddots
\end{array}\right),
$$
for all $t\in 1+p{\Bbb Z}_p$, where $t$ appears in the $i$-th position. Multiplying $\gamma'=(\gamma'_{ij})$
by suitable elements $h_{{i,i+1}}(t)$ with $t\in 1+p{\Bbb Z}_p$, we can transform $\gamma'$ into a 
matrix $\gamma''=(\gamma''_{ij})$, whose diagonal entries all equal $1$, i.e. $\gamma''$ is an upper unipotent matrix. Using left multiplication by elements $x_{\alpha_{ij}}(t)$ 
with $i<j$ and $t\in{\Bbb Z}_p$ suitably chosen, we see that $\gamma''$ can be transformed into the unit matrix. Since the transforming elements $x_{{\alpha}_{ij}}(t)$ and $h_{i+1,i}(t)$, which we used, are all contained in $K_*(p)$, this finally 
shows that $\gamma$ is contained in $K_*(p)$. Hence, the inclusion "$\supseteq$" holds and the claim is proven.

\bigskip

{\bf (4.3)  Irreducible representations of split semi-simple groups. } We fix a semi-simple algebraic ${\Bbb Q}$-group ${\bf G}$, which is split over 
${\Bbb Q}$. Hence, ${\bf G}={\bf G}_{{\lambda_0}}$ for some finite dimensional, irreducible representation $(\rho_{\lambda_0},L_{\lambda_0})$ of the Lie algebra 
${\mathfrak g}$ of ${\bf G}$. We recall the description of the irreducible representations of ${\bf G}$. If $K/{\Bbb Q}$ is an arbitrary extension field 
we define for $t\in K^*$ and $\alpha\in \Delta$ the elements $h_\alpha(t)$ as in [B], 3.2 (1), p. 13. The 
algebraic group ${\bf G}$ contains a maximal torus ${\bf T}$ such that for every extension $K/{\Bbb Q}$ the group of $K$-rational points ${\bf T}(K)$ is generated by the elements 
$$
h_\alpha(t), \quad t\in K^*,\,\alpha\in\Delta
$$
(cf. [B] 3.2 (1), p. 13). Moreover, ${\bf T}$ is defined over ${\Bbb Q}$ and splits over 
${\Bbb Q}$ (cf. [B], 3.3 (3), p. 15). To any $\lambda\in \Gamma_{\lambda_0}\,(=\langle \Pi_{\lambda_0} \rangle_{\Bbb Z})$ we attach a character $\lambda^\circ\in X({\bf T}(\bar{\Bbb Q}))$ by setting
$$
\lambda^\circ(\prod_{\alpha\in \Delta} h_\alpha(t_\alpha))=\prod_{\alpha\in \Delta} t_\alpha^{\lambda(h_\alpha)}\leqno(5)
$$
for all $t_\alpha\in \bar{\Bbb Q}^*$. 
In particular, $\lambda^\circ$ defines an algebraic character of ${\bf T}$ and the assignment 
$$
\begin{array}{cccl}
\varphi:&\Gamma_{\lambda_0}& \rightarrow& {\bf Mor}_{{\Bbb Q}-{\rm alg}}({\bf T},{\Bbb G}_m)=X({\bf T})\\
&\lambda&\mapsto&\lambda^\circ\\
\end{array}\leqno(6)
$$
is an isomorphism (cf. [B], 3.3 (3), p. 15). We note that $\mu^\circ$ is in fact the "exponential" of the weight $\mu\in \Pi_{\lambda_0}$: for $h=\prod_{\alpha\in \Delta} h_\alpha(t_\alpha)\in{\bf T}(\bar{\Bbb Q})$ and
any $v_\mu\in L_{\lambda_0}(\mu)$ it holds that
$$
h v_\mu=\mu^\circ(h) v_\mu\leqno(7)
$$
(cf. [B] 3.2 (1), p. 13). We call a weight $\lambda\in \Gamma_{\rm sc}$ analytically integral (for ${\bf G}$) precisely if $\Gamma_\lambda\subseteq\Gamma_{\lambda_0}$. Thus, if $\lambda$ 
is analytically integral then $\lambda^\circ\in X({\bf T})$ is defined. We call an algebraic weight $\lambda^\circ\in X({\bf T})$ dominant if the corresponding weight 
$\lambda\in \Gamma_{\lambda_0}$ is dominant.

\medskip

Let $\lambda\in\Gamma_{\rm sc}$ be an analytically integral weight.
%
%
For any root $\alpha\in \Phi$ and any $t\in \bar{\Bbb Q}$ we set
$$
x_\alpha^{\rho_\lambda}(t)=\exp t\rho_{\lambda}(x_\alpha) \in {\rm SL}(L_\lambda(\bar{\Bbb Q})).\leqno(8)
$$
There is a surjective morphism of groups 
$$
\pi_{\lambda^\circ,\bar{\Bbb Q}}:\,G_{\lambda_0,\bar{\Bbb Q}}\rightarrow G_{\lambda,\bar{\Bbb Q}}\subseteq {\rm SL}(L_\lambda(\bar{\Bbb Q})), 
$$
which sends $x_\alpha(t)$ to $x_\alpha^{\rho_\lambda}(t)$ for all $\alpha\in \Phi$ and $t\in \bar{\Bbb Q}$ (cf. [B], 3.2 (4), p. 14). Hence, we obtain a representation
$$
\pi_{\lambda^\circ,\bar{\Bbb Q}}:\,G_{\lambda_0,\bar{\Bbb Q}}\rightarrow {\rm SL}(L_\lambda(\bar{\Bbb Q})).\leqno(9)
$$
We recall from (4.1) that ${\bf G}_{\lambda_0}$ has a natural ${\Bbb Z}$-structure. 
Following [B], section  4.3, p. 22, equation (9) defines a representation of the algebraic group ${\bf G}={\bf G}_{\lambda_0}$, i.e. there is a morphism 
of algebraic groups
$$
\pi_{\lambda^\circ}:\,{\bf G}/{\Bbb Z}\rightarrow {\bf SL}(L_\lambda({\Bbb Z}))={\bf SL}_m/{\Bbb Z},\leqno(10)
$$ 
which is defined over ${\Bbb Z}$ with respect to the natural ${\Bbb Z}$-structure on ${\bf G}_{\lambda_0}$ and which on $\bar{\Bbb Q}$-points is given by 
$$
\pi_{\lambda^\circ}(x_\alpha(t))=x_\alpha^{\rho_\lambda}(t)\in {\rm SL}(L_\lambda(\bar{\Bbb Q}))\qquad (\alpha\in\Phi,t\in \bar{\Bbb Q}).\leqno(11)
$$
$(\pi_{\lambda^\circ},L_\lambda)$ is the irreducible algebraic representation of ${\bf G}$ of highest weight $\lambda^\circ\in X({\bf T})$ and any irreducible 
representation of ${\bf G}$ is
isomorphic to some $\pi_{\lambda^\circ}$ (note that if $(\pi_\lambda,L_\lambda)$ is the irreducible representation of highest weight $\lambda$ of ${\bf G}_{\lambda_0}$ then necessarily $\Gamma_{\lambda_0}\supseteq \Gamma_\lambda$).
We note that equation (10) implies that $\pi_{\lambda^\circ}$ defines a representation on $R$-points
$$
\pi_{\lambda^\circ}:\,{\bf G}(R)\rightarrow {\rm SL}(L_\lambda(R))={\rm SL}_m(R)\leqno(12)
$$
for any ${\Bbb Z}$-algebra $R$. In particular, $G_{\lambda_0,R}$ leaves $L_\lambda(R)$ invariant, 
i.e. $L_\lambda(R)$ is a $G_{\lambda_0,R}$-module. 


\bigskip

{\bf (4.4) The truncation of an irreducible representation of a split, semi-simple group. } As in section (4.3) we let ${\bf G}$ be an arbitrary ${\Bbb Q}$-split, semi-simple 
algebraic group, hence, ${\bf G}={\bf G}_{\lambda_0}$ for some dominant and integral weight $\lambda_0\in \Gamma_{\rm sc}$. Moreover, we fix a prime $p\in{\Bbb N}$ and we let 
$\lambda\in\Gamma_{\rm sc}$ be an analytically integral and dominant weight, i.e. $\Gamma_{\lambda_0}\supseteq \Gamma_\lambda$. Equation (12) then yields a representation on 
${\Bbb Z}_p$-points
$$
\pi_{\lambda^\circ}:\,{\bf G}({\Bbb Z}_p)\rightarrow{\rm SL}(L_\lambda({\Bbb Z}_p)).
$$
We look more closely at the quotient $L_\lambda({\Bbb Z}_p)/L_\lambda({\Bbb Z}_p,r)$. To this end, we set ${\cal S}_p={\Bbb Z}_p\otimes{\cal S}$. 
$L_\lambda({\Bbb Z}_p)\,(={\Bbb Z}_p\otimes L_\lambda({\Bbb Z}))$ is a ${\cal S}$-module with ${\cal S}$ acting on the second component of the tensor 
product. Moreover, ${\Bbb Z}_p$ acts on $L_\lambda({\Bbb Z}_p)$ via the first component of the tensor product 
and since the two actions commute we see that $L_\lambda({\Bbb Z}_p)$ is a ${\cal S}_p$-module.
(1.4) Lemma implies that $L_\lambda({\Bbb Z}_p,r)\,(={\Bbb Z}_p\otimes L_\lambda({\Bbb Z},r))$ is a ${\cal S}_p$-invariant submodule. In particular, the truncation
$$
{\bf L}^{[r]}_\lambda({\Bbb Z}_p)=L_\lambda({\Bbb Z}_p)/L_\lambda({\Bbb Z}_p,r)
$$ 
is a ${\cal S}_p$-module. We note that  $L_\lambda({\Bbb Z}_p)/L_\lambda({\Bbb Z}_p,r)\cong {\Bbb Z}_p\otimes(L_\lambda({\Bbb Z})/L_\lambda({\Bbb Z},r))$ as ${\Bbb Z}_p$-modules 
under the map $f$ which sends the class of $z\otimes v$ in $L_\lambda({\Bbb Z}_p)$ to $z\otimes (v+ L_\lambda({\Bbb Z},r))$. 
As above ${\cal S}_p$ acts componentwise on the factors of the tensor product ${\Bbb Z}_p\otimes(L_\lambda({\Bbb Z})/L_\lambda({\Bbb Z},r))$ and $f$ obviously 
is ${\cal S}_p$-equivariant. Thus, the isomorphism in (3.3) Proposition extends (under the respective assumptions) to an isomorphism of ${\cal S}_p$-modules
$$
{\rm id}\otimes\Phi:\,  \frac{L_\lambda({\Bbb Z}_p)}{L_\lambda({\Bbb Z}_p,r)}  \rightarrow \frac{L_{\lambda'}({\Bbb Z}_p)}{L_{\lambda'}({\Bbb Z}_p,r)} \leqno(13)
$$
and we obtain: if $\lambda,\lambda'\in\Gamma_{\rm sc}$ satisfy the conditions of (3.3) Proposition then there is an isomorphism of ${\cal S}_p$-modules 
$$
{\bf L}^{[r]}_\lambda({\Bbb Z}_p)\cong{\bf L}^{[r]}_{\lambda'}({\Bbb Z}_p).\leqno(14)
$$
%
Let now $x_\alpha(t)$ be one of the generators of $K_*(p)$, i.e. $t\in{\Bbb Z}_p$ if $\alpha$ is negative and $t\in p^{{\rm ht}(\alpha)}{\Bbb Z}_p$ if $\alpha$ is positive. Since 
$$
\pi_{\lambda^\circ}(x_\alpha(t))=\sum_{n\ge 0} \rho_\lambda(\frac{t^n x_\alpha^n}{n!})\leqno(15)
$$
and since $t^n \frac{x_\alpha^n}{n!}$ is contained in $S_p$ we obtain that the submodule $L_\lambda({\Bbb Z}_p,r)$ in particular is invariant under $K_*(p)$;
hence, ${\bf L}^{[r]}_\lambda({\Bbb Z}_p)$ is a $K_*(p)$-module. Moreover, since equation (15) holds with $\lambda$ replaced by $\lambda'$ we obtain that ${\bf L}^{[r]}_{\lambda'}({\Bbb Z}_p)$ also is a $K_*(p)$-module and since the isomorphism (14) is ${\cal S}_p$-equivariant we see that it is even is $K_*(p)$-equivariant.
%
%
Thus, we have proven the following

\bigskip

{\bf Proposition. }{\it 1.) Let $\lambda\in\Gamma_{\rm sc}$ be an analytically integral and dominant weight. The submodule $L_\lambda({\Bbb Z}_p,r)$ of 
$L_\lambda({\Bbb Z}_p)$ is $K_*(p)$-invariant, hence, ${\bf L}^{[r]}_\lambda({\Bbb Z}_p)$ is a $K_*(p)$-module.

\medskip

2.) Let $r\in {\Bbb N}$, let ${\cal T}\subseteq\Delta$ be a subset and let $\lambda,\lambda'\in\Gamma_{\rm sc}$ be analytically integral and dominant weights. 
We set $m_\alpha=\lambda(h_\alpha)$, $m_\alpha'=\lambda'(h_\alpha)$, $\alpha\in \Delta$. If

\begin{itemize}

\item $m_\alpha=m_\alpha'$ for all $\alpha\in \Delta-{\cal T}$

\item $m_\alpha>r$, $m_\alpha'>r$ for all $\alpha\in {\cal T}$

\item $\lambda\equiv \lambda'\pmod{p^{\lceil\frac{p}{p-1}r\rceil}\Gamma_{\rm sc}}$ 

\end{itemize}

then 
$$
{\bf L}^{[r]}_\lambda({\Bbb Z}_p)\cong{\bf L}^{[r]}_{\lambda'}({\Bbb Z}_p)
$$
as $K_*(p)$-modules.

}

\bigskip

{\bf (4.5) Remark. } 1.) Let $\Gamma\le {\bf G}({\Bbb Z})$ be an arithmetic subgroup which satisfies the following local condition at $p$:
$$
\Gamma\le K_*(p).\leqno(16)
$$
(4.4) Proposition then in particular implies that 
$$
{\bf L}^{[r]}_\lambda({\Bbb Z}_p)\cong{\bf L}^{[r]}_{\lambda'}({\Bbb Z}_p)
$$ 
as $\Gamma$-modules. As before, we call the 
$\Gamma$-module ${\bf L}^{[r]}_\lambda({\Bbb Z}_p)$ the truncation of length $r$ of the $\Gamma$-module $L_\lambda({\Bbb Z}_p)$.

\medskip

2.) If $\Gamma\le {\bf G}({\Bbb Z})$ satisfies equation (16) then $L_\lambda({\Bbb Z}_p,r)$ is $\Gamma$-invariant and we obtain
$$
\Gamma L_\lambda({\Bbb Z},r)\le L_\lambda({\Bbb Z})\cap L_\lambda({\Bbb Z}_p,r)=L_\lambda({\Bbb Z},r).
$$
Hence, ${\bf L}^{[r]}_\lambda({\Bbb Z})=L_\lambda({\Bbb Z})/L_\lambda({\Bbb Z},r)$ is a $\Gamma$-module and (4.4) Proposition shows that there is a 
${\bf G}(\Bbb Z)\cap K_*(p)$-equivariant isomorphism
$$
{\bf L}^{[r]}_\lambda({\Bbb Z})\cong {\bf L}^{[r]}_{\lambda'}({\Bbb Z}).
$$ 
In particular, the isomorphism is $\Gamma$-equivariant, i.e. (4.4) Proposition holds over ${\Bbb Z}$.

\section{The Hecke operator on Cohomology}

We introduce the Hecke operator acting on cohomology groups of $\Gamma$ with coefficients in irreducible representations of 
algebraic groups. To this end, from now on we will consider reductive groups and their representations. Although we will later only need the results about the 
Hecke operator acting on cohomology of $\Gamma$ with local coefficients, i.e. with coefficients given by modules over ${\Bbb Q}_p$ or ${\Bbb Z}_p$, we shall 
describe the global case, i.e. cohomology with coefficients in modules over ${\Bbb Q}$ or ${\Bbb Z}$; this is (at least formally) a stronger statement. In (5.8)  
we explain how the corresponding local results are easy consequences of the global ones. Alternatively, one can deduce the local results by repeating the arguments 
from the global case.

\bigskip

{\bf (5.1) Irreducible representations of Reductive groups and their truncations. }

\bigskip

{\bf (5.1.1) Reductive groups. } We let 
$\tilde{\bf G}/{\Bbb Q}$ be a connected reductive algebraic group, which splits over ${\Bbb Q}$. We denote by $\tilde{\bf T}$ a maximal ${\Bbb Q}$-split 
torus in $\tilde{\bf G}$ and we denote by ${\bf G}=(\tilde{\bf G},\tilde{\bf G})$ the derived group. ${\bf G}$ is a semi simple ${\Bbb Q}$-split group 
and $\tilde{\bf T}$ contains a maximal ${\Bbb Q}$-split torus ${\bf T}$ of ${\bf G}$ (cf. [S], Prop. 8.1.8 (iii), p. 135).
 We denote by ${\mathfrak g}$ the Lie algebra of ${\bf G}$, hence, ${\bf G}={\bf G}_{\rho_0}$ for some faithful representation $(\rho_{\lambda_0},L_{\lambda_0})$ of ${\mathfrak g}$ of highest weight $\lambda_0\in \Gamma_{\rm sc}$. We choose the Cartan subalgebra ${\mathfrak h}$ in ${\mathfrak g}$ and the basis 
$\Delta$ of the root system of $({\mathfrak g},{\mathfrak h})$ such that ${\bf T}$ is given as in (4.3). 
We know that
$$
\tilde{\bf G}(\bar{\Bbb Q})={\bf G}(\bar{\Bbb Q})\times_{\sf g} {\Bbb G}_m^a(\bar{\Bbb Q}),
$$
where ${\Bbb G}_m^a(\bar{\Bbb Q})={\rm Rad}(\tilde{\bf G}(\bar{\Bbb Q}))$ is the radical of $\tilde{\bf G}(\bar{\Bbb Q})$ and ${\sf g}\le {\bf G}(\bar{\Bbb Q})\cap {\Bbb G}_m^a(\bar{\Bbb Q})$ is a finite subgroup. We note that ${\Bbb G}_m^a(\bar{\Bbb Q})$ is contained in the centre of 
$\tilde{\bf G}(\bar{\Bbb Q})$, hence, ${\sf g}$ is contained in the centre of ${\bf G}(\bar{\Bbb Q})$ (cf. [S], Prop. 7.3.1, p. 120 and Cor. 8.1.6, p. 
134). In particular, we obtain
$$
\tilde{\bf T}(\bar{\Bbb Q})={\bf T}(\bar{\Bbb Q})\times_{\sf g} {\Bbb G}_m^a(\bar{\Bbb Q}),
$$
i.e. any $t\in\tilde{\bf T}(\bar{\Bbb Q})$ can be written
$$
t=t^0z,\leqno(1)
$$
where $t^0\in {\bf T}(\bar{\Bbb Q})$ and $z\in {\Bbb G}_m^a(\bar{\Bbb Q})$. The group of algebraic ${\Bbb Q}$-characters of $\tilde{\bf T}$ is given as
$$
X(\tilde{\bf T})=\{\lambda^\circ\otimes\kappa,\,\lambda^\circ\in X({\bf T}), \kappa\in X({\Bbb G}_m^a):\;\lambda^\circ|_{\sf g}=\kappa|_{\sf g}\};
$$
explicitly, the ${\Bbb Q}$-character $\lambda^\circ\otimes\kappa$ is defined on $\bar{\Bbb Q}$-points as 
$$
\lambda^\circ\otimes\kappa(t^0 z):=\lambda^\circ(t^0)\kappa(z),\leqno(2)
$$ 
where $t^0\in {\bf T}(\bar{\Bbb Q})$ and $z\in {\Bbb G}_m^a(\bar{\Bbb Q})$. In particular, since the simple roots $\alpha^\circ\in X({\bf T})$, 
$\alpha\in\Delta$ (cf. equation (5) in section (4.3)), vanish on the centre of ${\bf G}(\bar{\Bbb Q})$, hence, on ${\sf g}$, they extend to characters of 
$\tilde{\bf T}(\bar{\Bbb Q})$ by defining $\alpha^\circ(t^0 z)=\alpha^\circ(t^0)$, 
$t^0\in {\bf T}(\bar{\Bbb Q})$, $z\in {\Bbb G}_m^a(\bar{\Bbb Q})$. 

\medskip

We note that the characters $\lambda^\circ$, $\lambda\in \Gamma_{\lambda_0}$, are defined over ${\Bbb Z}$ by equation (5) in section (4.3), hence, the
elements in $X(\tilde{\bf T})$ are defined over ${\Bbb Z}$. 


\bigskip

{\bf (5.1.2) Irreducible representations. } Let $\tilde{\lambda}\in X(\tilde{\bf T})$ be an algebraic weight. Thus, $\tilde{\lambda}$ decomposes
$$
\tilde{\lambda}=\lambda^\circ\otimes\kappa, \leqno(3)
$$ 
where $\kappa=\tilde{\lambda}|_{{\Bbb G}_m^a}$ and $\lambda^\circ=\tilde{\lambda}|_{{\bf T}}\in X({\bf T})$ is the image of a weight $\lambda\in\Gamma_{\lambda_0}$ under the map 
$\lambda\mapsto\lambda^\circ$ (cf. equation (6) in section (4.3)). We call $\tilde{\lambda}$ dominant if $\lambda$ is dominant.

\medskip

Let $\tilde{\lambda}\in X(\tilde{\bf T})$ be dominant. We denote by
$$
\pi_{\tilde{\lambda}}:\,\tilde{\bf G}\rightarrow {\bf GL}(L_{\tilde{\lambda}}({\Bbb Z}))={\bf GL}_m
$$ 
the irreducible representation of the algebraic group $\tilde{\bf G}$ of highest weight $\tilde{\lambda}$ (as in section (4.3) we identify the algebraic groups ${\bf GL}(L_{\tilde{\lambda}}({\Bbb Z}))={\bf GL}_m/{\Bbb Z}$, where $m={\rm dim}\,L_{\tilde{\lambda}}({\Bbb Z})$). The representation $(\pi_{\tilde{\lambda}},L_{\tilde{\lambda}})$ is defined over ${\Bbb Z}$; it is given as 
$$
\pi_{\tilde{\lambda}}=\pi_{\lambda^\circ}\otimes \kappa,
$$
where $\lambda^\circ$ and $\kappa$ are as in equation (3) and $\pi_{\lambda^\circ}:\,{\bf G}/{\Bbb Z}\rightarrow {\bf SL}(L_{\lambda}({\Bbb Z}))$ is the 
irreducible representation of ${\bf G}$ of highest weight $\lambda^\circ$ (cf. equation (10) in section (4.3); note that $\lambda\in\Gamma_{\lambda_0}$, 
hence, $\Gamma_{\lambda_0}\supseteq \Gamma_\lambda$).  Explicitly, this means
$$
\pi_{\tilde{\lambda}}(g)=\kappa(z)\,\pi_{\lambda^\circ}(g^0)\leqno(4)
$$
for any $g=g^0z\in\tilde{\bf G}(\bar{\Bbb Q})$ with $g^0\in {\bf G}(\bar{\Bbb Q})$ and $z\in {\Bbb G}_m^a(\bar{\Bbb Q})$. The representation $\pi_{\tilde{\lambda}}$ is well defined because $\lambda^\circ$ and $\kappa$ coincide on ${\sf g}$. 
%
%

\bigskip

{\bf (5.1.3) Truncations of irreducible representations of reductive groups. } Since $\kappa$ is a character, equation (4) implies that the representation spaces 
$L_{\tilde{\lambda}}$ and $L_\lambda$ of $\pi_{\tilde{\lambda}}$ and $\pi_{\lambda^\circ}$ coincide, i.e. for any ${\Bbb Z}$-algebra $R$ we have an 
equality 
\begin{itemize}

\item $L_{\tilde{\lambda}}(R)=L_\lambda(R)$

\end{itemize} 

of $R$-modules and also of ${\bf G}(R)$-modules. In contrast to $L_\lambda$, on $L_{\tilde{\lambda}}$ in addition we have an action of ${\Bbb G}_m^a$ via 
the character $\kappa$ (cf. equation (3)).  
Any ${\bf T}(R)$ resp. ${\bf G}(R)$-invariant subspace of $L_\lambda(R)$ therefore also 
is a subspace of $L_{\tilde{\lambda}}(R)$ and this subspace is $\tilde{\bf T}(R)$ resp. $\tilde{\bf G}(R)$-invariant. In particular, we define the 
following subspaces of $L_{\tilde{\lambda}}(R)$

\begin{itemize}

\item $L_{\tilde{\lambda}}(R,\mu):=L_\lambda(R,\mu)$, $\mu\in\Pi_\lambda$

\item $L_{\tilde{\lambda}}(R,r):=L_\lambda(R,r)$, $r\in{\Bbb N}_0$.

\end{itemize}

Thus, $L_{\tilde{\lambda}}(R,\mu)$ is invariant under $\tilde{\bf T}(R)$ and it is the weight subspace of $L_{\tilde{\lambda}}(R)$ of weight 
$$
\mu^\circ\otimes\kappa\leqno(5)
$$ 
with respect to $\tilde{\bf T}(R)$ (cf. equation (7) in section (4.3)) and $L_{\tilde{\lambda}}(R,r)$ is a $\Gamma$-invariant subspace of $L_{\tilde{\lambda}}(R)$. 
%
%
%
In particular, we obtain the weight decompositions with 
respect to the torus $\tilde{\bf T}$ 
$$
L_{\tilde{\lambda}}(R)=\bigoplus_{\mu\in\Pi_\lambda} L_{\tilde{\lambda}}(R,\mu).\leqno(6)
$$
and (using section (1.2))
$$
L_{\tilde{\lambda}}(R,r)=\bigoplus_{\mu\in\Pi_\lambda \atop 0\le {\rm ht}_\lambda(\mu)\le r} p^{r-{\rm ht}_\lambda(\mu)}\,L_{\tilde{\lambda}}(R,\mu) \oplus
\bigoplus_{\mu\in\Pi_\lambda \atop {\rm ht}_\lambda(\mu)> r} L_{\tilde{\lambda}}(R,\mu).
$$
The quotient
$$
{\bf L}^{[r]}_{\tilde{\lambda}}({R})=\frac{L_{\tilde{\lambda}}({R})}{L_{\tilde{\lambda}}({R,}r)}
$$
is $\Gamma$-invariant and as before we call the $\Gamma$-module ${\bf L}^{[r]}_{\tilde{\lambda}}(R)$ the truncation of $L_{\tilde{\lambda}}(R)$ of 
length $r$. Again, we will need this only in the cases $R={\Bbb Z},{\Bbb Z}_p$.

\bigskip

{\it Remark. } Of course, since $\Gamma\le {\bf G}({\Bbb Z})$ we obviously have $L_{\tilde{\lambda}}({\Bbb Z}_p,r)=L_\lambda({\Bbb Z}_p,r)$ and 
$$
{\bf L}^{[r]}_{\tilde{\lambda}}({R})={\bf L}^{[r]}_\lambda({R})
$$ 
as $\Gamma$-modules. On the other hand, $L_{\tilde{\lambda}}({\Bbb Z})$ and $L_{\tilde{\lambda}}({\Bbb Z},r)$ are submodule of 
$L_{\tilde{\lambda}}({\Bbb Q})$, hence, the torus $\tilde{\bf T}({\Bbb Q})$ acts on them. This will yield Hecke operators ${\Bbb T}(h)$ for 
elements $h\in\tilde{\bf T}({\Bbb Q})$ which act on the cohomology of $\Gamma$ with coefficients in $L_{\tilde{\lambda}}({\Bbb Z}_p)$ and 
$L_{\tilde{\lambda}}({\Bbb Z}_p,r)$ and, hence, in ${\bf L}^{[r]}_{\tilde{\lambda}}({\Bbb Z}_p)$ (cf. below). In contrast, on 
$H^\bullet(\Gamma,{\bf L}^{[r]}_\lambda({\Bbb Z}_p))$ there is no action of ${\Bbb T}(h)$ for elements $h\in\tilde{\bf T}({\Bbb Q})$.


\bigskip

{\bf (5.2) The Hecke operator acting on cohomology. } In the remainder of section 5, we fix a ${\Bbb Q}$-split reductive group $\tilde{\bf G}$. We use 
the notations introduced in section (5.1), e.g. ${\bf G}=(\tilde{\bf G},\tilde{\bf G})$ is the derived group of $\tilde{\bf G}$ and ${\bf G}={\bf G}_{\lambda_0}$ 
for some irreducible representation $(\rho_{\lambda_0},L_{\lambda_0})$ of the Lie algebra ${\mathfrak g}$ of ${\bf G}$. We recall that 
${\bf G}$ has a natural ${\Bbb Z}$-structure ${\bf G}/{\Bbb Z}={\bf G}_{\lambda_0}/{\Bbb Z}$ (cf. (4.1)). Moreover, we fix a prime $p\in{\Bbb N}$ and 
an arithmetic subgroup $\Gamma\le {\bf G}({\Bbb Z})$ such that $\Gamma$ satisfies the following local condition at $p$:
$$
\Gamma\le K_*(p).
$$
We denote by $\tilde{\lambda}\in X(\tilde{\bf T})$ a dominant algebraic weight and we write $\tilde{\lambda}$ as 
$$
\tilde{\lambda}=\lambda^\circ\otimes\kappa
$$ 
where $\lambda^\circ$ and $\kappa$ are as in equation (3). The highest weight module $(\pi_{\tilde{\lambda}},L_{\tilde{\lambda}}({\Bbb Z}))$ of $\tilde{\bf G}$
decomposes $\pi_{\tilde{\lambda}}=\pi_{\lambda^\circ}\otimes\kappa$ (cf. equation (4)). Since $\pi_{\tilde{\lambda}}$ is defined over ${\Bbb Z}$, the group 
$\Gamma\le {\bf G}({\Bbb Z})$ leaves the lattice $L_{\tilde{\lambda}}({\Bbb Z})=L_\lambda({\Bbb Z})$ invariant, i.e. $L_{\tilde{\lambda}}(R)$ is a 
$\Gamma$-module for any ${\Bbb Z}$-algebra $R$ and the cohomology groups
$$
H^i(\Gamma,L_{\tilde{\lambda}}(R))
$$
are defined.

\medskip

The $\tilde{\bf T}({\Bbb Q})$-module structure on $L_{\tilde{\lambda}}({\Bbb Q})$ enables us to
define for elements $h\in \tilde{\bf T}({\Bbb Q})$ a Hecke operator acting on $H^i(\Gamma,L_{\tilde{\lambda}}(\Bbb Q))$.
%
%
To explain this, we fix {\it strictly} positive integers $e_\alpha\in{\Bbb N}$, $\alpha\in \Delta$, 
and an element in $h\in\tilde{\bf T}({\Bbb Q})$, which satisfies 
$$
\alpha^\circ(h)=p^{e_\alpha}\leqno(7)
$$
for all $\alpha\in\Delta$. 
%
%
We set
$$
T(h)=\Gamma h\Gamma. 
$$
The double coset $T(h)$ acts on the cohomolgy of $\Gamma$ with values in $L_{\tilde{\lambda}}({\Bbb Q})$ (more generally, $T(h)$ acts on the cohomology 
of any $\langle \Gamma,h^{-1}\rangle_{\rm semi}$-module $V$ where $\langle \Gamma,h^{-1}\rangle_{\rm semi}\le \tilde{\bf G}({\Bbb Q})$ is the sub semigroup generated by 
$\Gamma$ and $h^{-1}$). 
To recall the definition of this action we fix a system of representatives $\gamma_1,\ldots,\gamma_d$ for $h^{-1}\Gamma h\cap \Gamma\backslash \Gamma$, hence, we obtain
$$
T(h)=\bigcup_{i=1,\ldots,d} \Gamma h \gamma_i.
$$
For any $\eta\in\Gamma$ and any index $i$ satisfying $1\le i\le d$ there is an index $\eta(i)$ such that 
$$
\Gamma h \gamma_i\eta=\Gamma h \gamma_{\eta(i)},
$$ 
In particular, there are $\rho_i(\eta)\in \Gamma$, $i=1,\ldots,d$, such that 
$h \gamma_i\eta=\rho_i(\eta) h \gamma_{\eta(i)}$. Let now $c\in C^h(\Gamma,L_{\tilde{\lambda}}({\Bbb Q}))$ be any cochain; we then define 
$T(h)(c)$ as the cochain $c'\in C^h(\Gamma,L_{\tilde{\lambda}}({\Bbb Q}))$, which is given by 
$$
c'(\eta_0,\ldots,\eta_h)=\sum_{1\le i\le d} (h\gamma_i)^{-1} c(\rho_i(\eta_0),\ldots,\rho_i(\eta_h)\leqno(8)
$$
(cf. [K-P-S], p. 227). Thus, $T(h)$ defines an operator on $C^\bullet(\Gamma,L_{\tilde{\lambda}}({\Bbb Q}))$. Since the action of $T(h)$ commutes with the coboundary 
operator, $T(h)$ acts on cohomology with coefficients in $L_{\tilde{\lambda}}({\Bbb Q})$, i.e. $T(h)$ defines an element in ${\rm End}(H^\bullet(\Gamma,L_{\tilde{\lambda}}({\Bbb Q})))$. 


\bigskip

{\it Note. } The Hecke algebra acts from the right on cohomology 
but since we only consider a single Hecke operator $T$ this will not become relevant and we therefore write the image of a cohomology class $c$ under $T$ as 
$T(c)$.

\bigskip

{\it Remark. } The action of $T(h)$ on cochains $C^\bullet(\Gamma,L_{\tilde{\lambda}}({\Bbb Q}))$ may still depend on the choice of the system of representatives 
$\gamma_1,\ldots,\gamma_d$ but the operator $T(h)$ on $H ^\bullet(\Gamma,L_{\tilde{\lambda}}({\Bbb Q}))$ is independent of the choice of these representatives (cf. [K-P-S], p. 227).

\medskip


\bigskip

{\bf (5.3) Normalization of Hecke operators acting on cohomology. } In general, the action of $T(h)$ on cochains does not leave the lattice 
$C^h(\Gamma,L_{\tilde{\lambda}}({\Bbb Z}))$ in  $C^h(\Gamma,L_{\tilde{\lambda}}({\Bbb Q}))$ invariant, i.e. $T(h)$ does not act on the subcomplex 
$C^\bullet(\Gamma,L_{\tilde{\lambda}}({\Bbb Z}))$ of $C^\bullet(\Gamma,L_{\tilde{\lambda}}({\Bbb Q}))$. To achieve 
this, we have to normalize the Hecke operator as follows. For any dominant $\tilde{\lambda}\in X(\tilde{\bf T})$ we define the normalized 
operator ${\Bbb T}(h)$ as
$$
{\Bbb T}(h)=\tilde{\lambda}(h)\;T(h)\,\in {\rm End}(C^\bullet(\Gamma,L_{\tilde{\lambda}}({\Bbb Q}))).\leqno(9)
$$
We note that in the rank-$1$ case the normalized Hecke operator ${\Bbb T}$ corresponds to the classical Hecke operator on modular forms.

\medskip

\medskip

We want to study the normalized Hecke operator ${\Bbb T}(h)$ on cohomolgy with integral coefficients. This is based on the following

\bigskip

{\bf (5.4) Lemma. }{\it 1.) For any weight $\mu\in\Pi_\lambda$ we have 
$$
\tilde{\lambda}(h) \pi_{\tilde{\lambda}}(h^{-1}) L_{\tilde{\lambda}}({\Bbb Z},\mu)\subseteq p^{{\rm ht}_\lambda(\mu)} L_{\tilde{\lambda}}({\Bbb Z},\mu).
$$

2.) For any $\gamma\in \Gamma$ we have 
$$
\tilde{\lambda}(h) \pi_{\tilde{\lambda}}(h \gamma)^{-1} L_{\tilde{\lambda}}({\Bbb Z})\subseteq L_{\tilde{\lambda}}({\Bbb Z}).
$$

3.) For any $\gamma\in\Gamma$ we have 
$$
\tilde{\lambda}(h) \pi_{\tilde{\lambda}}(h \gamma)^{-1} L_{\tilde{\lambda}}({\Bbb Z},r)\subseteq p^r L_{\tilde{\lambda}}({\Bbb Z}).
$$

}

\medskip

{\it Proof. } %
1.) We write $\tilde{\lambda}=\lambda^\circ\otimes\kappa$ as in equation (3) and we write the element defining the Hecke operator as $h=h^0z$, $h^0\in {\bf T}(\bar{\Bbb Q})$, $z\in{\Bbb G}_m^a(\bar{\Bbb Q})$. Equation (7) implies that 
$$
\alpha^\circ(h^0)=p^{e_\alpha}\leqno(10)
$$ 
for all simple roots $\alpha\in\Delta$. 
We let $v_\mu\in L_{\tilde{\lambda}}({\Bbb Z},\mu)$ be arbitrary. Equation (5) implies that $v_\mu$ has weight $\mu^\circ\otimes\kappa$, hence, we obtain 
$$
\pi_{\tilde{\lambda}}(h)v_\mu=\kappa(z) \mu^\circ(h^0) v_\mu.
$$
Since $\mu\in\Pi_{\lambda}$, it has the form $\mu=\lambda-\sum_{\alpha\in \Delta} c_\alpha\alpha$ with all
$c_\alpha\in{\Bbb N}_0$ and we further obtain using equation (10)
\begin{eqnarray*}
\kappa(z) \mu^\circ(h^0)&=&\kappa(z) \lambda^\circ(h^0) \prod_{\alpha\in \Delta}
(\alpha^\circ)^{-c_\alpha}(h^0)\\
&=&\tilde{\lambda}(h) \prod_{\alpha\in \Delta} (\alpha^\circ)^{-c_\alpha}(h^0)\\
&=&\tilde{\lambda}(h) \prod_{\alpha\in \Delta} p^{-e_\alpha c_\alpha}.\\
\end{eqnarray*}
Thus, we obtain
$$
\tilde{\lambda}(h)\pi_{\tilde{\lambda}}(h^{-1}) v_\mu=\prod_{\alpha\in \Delta} p^{e_\alpha c_\alpha} v_\mu.
$$
Since $e_\alpha$ is strictly positive (cf. equation (7)) we deduce that 
$$
\sum_{\alpha\in\Delta} c_\alpha e_\alpha\ge \sum_{\alpha\in\Delta} c_\alpha={\rm ht}_\lambda(\mu).
$$
Hence, we obtain $\tilde{\lambda}(h) \pi_{\tilde{\lambda}}(h^{-1}) v_\mu=c_\mu v_\mu$ with $c_\mu\in p^{{\rm ht}_\lambda(\mu)}{\Bbb Z}$.
This proves 1.) 

\medskip

2.) Since ${\rm ht}_\lambda(\mu)\ge 0$ for all $\mu\in\Pi_\lambda$, part 1.) implies that $\tilde{\lambda}(h) \pi_{\tilde{\lambda}}(h^{-1}) L_{\tilde{\lambda}}({\Bbb Z},\mu)
\subseteq L_{\tilde{\lambda}}({\Bbb Z},\mu)$ for all $\mu\in\Pi_\lambda$, hence, $\tilde{\lambda}(h) \pi_{\tilde{\lambda}}(h^{-1}) L_{\tilde{\lambda}}({\Bbb Z})
\subseteq L_{\tilde{\lambda}}({\Bbb Z})$ (cf. equation (6) in (5.1)). Since $\pi_{\tilde{\lambda}}(\gamma^{-1})$ leaves $L_{\tilde{\lambda}}({\Bbb Z})$ invariant, the claim follows.

\medskip

3.) Again, part 1.) and the weight decomposition of $L_{\tilde{\lambda}}({\Bbb Z},r)$ (cf. section (5.1.3)) imply that 
$$
\tilde{\lambda}(h) \pi_{\tilde{\lambda}}(h^{-1}) L_{\tilde{\lambda}}({\Bbb Z},r)\subseteq p^r L_{\tilde{\lambda}}({\Bbb Z}).
$$ 
Since $\pi_{\tilde{\lambda}}(\gamma^{-1})$ leaves $L_{\tilde{\lambda}}({\Bbb Z})$ invariant we obtain the claim. Thus, the proof of the Lemma is complete.

\bigskip

{\bf (5.5) Proposition. }{\it Let $\tilde{\lambda}\in X(\tilde{\bf T})$ be a dominant algebraic weight. The Hecke operator ${\Bbb T}(h)$ leaves the subcomplex $C^\bullet(\Gamma,L_{\tilde{\lambda}}({\Bbb Z}))$ of $C^\bullet(\Gamma,L_{\tilde{\lambda}}({\Bbb Q}))$ and, hence, the complex
$C^\bullet(\Gamma,L_{\tilde{\lambda}}({\Bbb Z}/(p^r)))=C^\bullet(\Gamma,L_{\tilde{\lambda}}({\Bbb Z}))\otimes{\Bbb Z}/(p^r)$ invariant.

}

\medskip

{\it Proof. }  We let $c\in C^h(\Gamma,L_{\tilde{\lambda}}({\Bbb Z}))$ be arbitrary and we decompose $\Gamma h\Gamma=\bigcup_i\Gamma h\gamma_i$. Hence, for any tuple $(\eta_0,\ldots,\eta_h)\in \Gamma^{h+1}$ we have
$$
c(\rho_i(\eta_0),\ldots,\rho_i(\eta_h))\in L_{\tilde{\lambda}}({\Bbb Z}).
$$
Thus, the definition of ${\Bbb T}(h)(c)$ in equation (8) together with (5.4) Lemma 2.) implies that 
$$
({\Bbb T}(h)(c)) (\eta_0,\ldots,\eta_h)\in\, L_{\tilde{\lambda}}({\Bbb Z}).
$$
Since $(\eta_0,\ldots,\eta_h)$ was arbitrary this implies that ${\Bbb T}(h)(c)$ is contained in $C^h(\Gamma,L_{\tilde{\lambda}}({\Bbb Z}))$, which proves the first claim. 
Since the Hecke operator ${\Bbb T}(h)$ leaves the subspace $C^\bullet(\Gamma,p^r L_{\tilde{\lambda}}({\Bbb Z}))=p^rC^\bullet(\Gamma,L_{\tilde{\lambda}}({\Bbb Z}))$ in $C^\bullet(\Gamma,L_{\tilde{\lambda}}({\Bbb Z}))$ invariant, ${\Bbb T}(h)$ also acts on 
$$
C^\bullet(\Gamma,L_{\tilde{\lambda}}({\Bbb Z})/p^r L_{\tilde{\lambda}}({\Bbb Z}))
=C^\bullet(\Gamma,L_{\tilde{\lambda}}({\Bbb Z}))/p^rC^\bullet(\Gamma,L_{\tilde{\lambda}}({\Bbb Z})),
$$
which implies the second claim. 
Thus, the proof of the Proposition is complete.

\bigskip

{\bf (5.6) Cohomology of the truncating module. } Let $R$ be a ${\Bbb Z}$-algebra. We denote by 
$L_{\tilde{\lambda}}({R}/(p^r),r)$ the image of the truncating module $L_{\tilde{\lambda}}({R},r)=L_{\lambda}({R},r)\le L_{\tilde{\lambda}}({R})$ under 
the canonical map given by mod $p^r$-reduction
$$
L_{\tilde{\lambda}}({R})\rightarrow L_{\tilde{\lambda}}(\frac{R}{(p^r)})
=L_{\tilde{\lambda}}({R})\otimes_{R}\frac{R}{(p^r)}\leqno(11)
$$
%
%
sending $v\mapsto v\otimes 1=v+p^r L_{\tilde{\lambda}}({R})$. Thus, $L_{\tilde{\lambda}}({R}/(p^r),r)$ is a submodule of $L_{\tilde{\lambda}}({R}/(p^r))$ and 
$$
L_{\tilde{\lambda}}(\frac{{R}}{(p^r)},r)=
\bigoplus_{\mu\in\Pi_\lambda  \atop  {\rm ht}_\lambda(\mu)\le r} p^{r-{\rm ht}_\lambda(\mu)}
L_{\tilde{\lambda}}({R},\mu) 
\oplus \bigoplus_{\mu\in\Pi_\lambda \atop {\rm ht}_\lambda(\mu) > r} L_{\tilde{\lambda}}({R},\mu)\pmod{p^r L_{\tilde{\lambda}}({R})}.
%
$$
%


\bigskip

{\bf Proposition. } {\it 1.) ${\Bbb T}(h)$ leaves the subcomplex 
$C^\bullet(\Gamma,L_{\tilde{\lambda}}({\Bbb Z},r))$ of $C^\bullet(\Gamma,L_{\tilde{\lambda}}({\Bbb Z}))$ and the subcomplex $C^\bullet(\Gamma,L_{\tilde{\lambda}}({\Bbb Z}/(p^r),r))$ of $C^\bullet(\Gamma,L_{\tilde{\lambda}}({\Bbb Z}/(p^r)))$ invariant. 

\medskip

2.) ${\Bbb T}(h)$ annihilates the subcomplex $C^\bullet(\Gamma,L_{\tilde{\lambda}}({\Bbb Z}/(p^r),r))$ of $C^\bullet(\Gamma,L_{\tilde{\lambda}}({\Bbb Z}/(p^r)))$.

}

\bigskip

We note that (5.5) Proposition implies that ${\Bbb T}(h)$ acts on $C^\bullet(\Gamma,L_{\tilde{\lambda}}({\Bbb Z}))$ and on $C^\bullet(\Gamma,L_{\tilde{\lambda}}({\Bbb Z}/(p^r)))$.  

\bigskip

{\it Proof. } We prove 1.) and 2.) together. Let $c\in C^h(\Gamma,L_{\tilde{\lambda}}({\Bbb Z},r))$ be arbitrary and write 
$\Gamma h\Gamma=\bigcup_i\Gamma h\gamma_i$. Hence, for any tuple $(\eta_0,\ldots,\eta_h)\in \Gamma^{h+1}$ we have
$$
c(\rho_i(\eta_0),\ldots,\rho_i(\eta_h))\in L_{\tilde{\lambda}}({\Bbb Z},r).
$$
The definition of ${\Bbb T}(h)(c)$ in equation (8) together with (5.4) Lemma 3.) therefore implies that
$$
{\Bbb T}(h)(c)(\eta_0,\ldots,\eta_h)\in p^r L_{\tilde{\lambda}}({\Bbb Z}).
$$
Since $p^r L_{\tilde{\lambda}}({\Bbb Z})\subseteq L_{\tilde{\lambda}}({\Bbb Z},r)$ this implies that ${\Bbb T}(h)(c)$ is contained in 
$C^h(\Gamma,L_{\tilde{\lambda}}({\Bbb Z},r))$, hence, ${\Bbb T}(h)$ leaves $C^h(\Gamma,L_{\tilde{\lambda}}({\Bbb Z},r))$ invariant. Moreover, 
since ${\Bbb T}(h)(c)$ has values in $p^r L_{\tilde{\lambda}}({\Bbb Z})$, it is contained in $C^\bullet(\Gamma,p^r L_{\tilde{\lambda}}({\Bbb Z}))$ and therefore 
vanishes in $C^\bullet(\Gamma,L_{\tilde{\lambda}}({\Bbb Z}/(p^r)))=C^\bullet(\Gamma,L_{\tilde{\lambda}}({\Bbb Z})/p^r L_{\tilde{\lambda}}({\Bbb Z}))$. 
Hence, ${\Bbb T}(h)$ annihilates the image of $C^\bullet(\Gamma,L_{\tilde{\lambda}}({\Bbb Z},r))$ in $C^\bullet(\Gamma,L_{\tilde{\lambda}}({\Bbb Z}/(p^r)))$. 
Thus, the proof is complete.

\bigskip

{\bf (5.7) } In the following diagram we display the relations between the several complexes that we used (induced by the corresponding 
relations between the coefficient systems)

$$
\begin{array}{ccccc}
&&C^\bullet(\Gamma,L_{\tilde{\lambda}}({\Bbb Q}))&&\\
&&\cup&&\\
C^\bullet(\Gamma,L_{\tilde{\lambda}}(\frac{\Bbb Z}{(p^r)}))&\leftarrow&C^\bullet(\Gamma,L_{\tilde{\lambda}}({\Bbb Z}))
&\rightarrow&C^\bullet(\Gamma,{\bf L}^{[r]}_{\tilde{\lambda}}({\Bbb Z}))\\
\cup&&\cup&&\\
C^\bullet(\Gamma,L_{\tilde{\lambda}}(\frac{\Bbb Z}{(p^r)},r))&\leftarrow&C^\bullet(\Gamma,L_{\tilde{\lambda}}({\Bbb Z},r)).&&\\
\end{array}
$$

\medskip

The Hecke operator ${\Bbb T}(h)$ acts on the complex $C^\bullet(\Gamma,L_{\tilde{\lambda}}({\Bbb Q}))$ and in (5.5) and (5.6) Propositions 
we have seen that ${\Bbb T}(h)$ leaves the subcomplexes $C^\bullet(\Gamma,L_{\tilde{\lambda}}({\Bbb Z}))$ and $C^\bullet(\Gamma,L_{\tilde{\lambda}}({\Bbb Z},r))$ invariant; 
in particular it also acts on $C^\bullet(\Gamma,L_{\tilde{\lambda}}({\Bbb Z}/(p^r)))$ and in (5.6) Proposition we have seen that ${\Bbb T}(h)$ leaves the subcomplex 
$C^\bullet(\Gamma,L_{\tilde{\lambda}}({\Bbb Z}/(p^r),r))$ of $C^\bullet(\Gamma,L_{\tilde{\lambda}}({\Bbb Z}/(p^r)))$ invariant. In particular, ${\Bbb T}(h)$ acts on the quotient
complex $C^\bullet(\Gamma,{\bf L}^{[r]}_{\tilde{\lambda}}({\Bbb Z}))\cong C^\bullet(\Gamma,L_{\tilde{\lambda}}({\Bbb Z}))/C^i(\Gamma,L_{\tilde{\lambda}}({\Bbb Z},r))$. 
Thus, ${\Bbb T}(h)$ acts on all of the above complexes and it annihilates $C^\bullet(\Gamma,L_{\tilde{\lambda}}({\Bbb Z}/(p^r),r))$. Hence, we have proven

\bigskip

{\bf Proposition. }{\it The normalized Hecke operator ${\Bbb T}(h)$ acts on the cohomology groups

\medskip

$H^\bullet(\Gamma,L_{\tilde{\lambda}}({\Bbb Z}))$, $H^\bullet(\Gamma,L_{\tilde{\lambda}}({\Bbb Z},r))$, $H^\bullet(\Gamma,L_{\tilde{\lambda}}(\frac{\Bbb Z}{(p^r)}))$, 
$H^\bullet(\Gamma,L_{\tilde{\lambda}}(\frac{{\Bbb Z}}{(p^r)},r))$ and $H^\bullet(\Gamma,{\bf L}^{[r]}_{\tilde{\lambda}}({\Bbb Z}))$  

\medskip

and it annihilates the second last one.
}

%
%
%
%
%
%
%
%
%
%
%

\bigskip

{\bf (5.8) Cohomology with local coefficients.  } Since $T(h)=\Gamma h\Gamma$ acts on $C^i(\Gamma,V)$ where $V$ is any $\langle \Gamma,h^{-1}\rangle_{\rm semi}$-module, $T(h)$
and, hence, ${\Bbb T}(h)$ acts on $C^i(\Gamma,L_{\tilde{\lambda}}({\Bbb Q}_p))=C^i(\Gamma,L_{\tilde{\lambda}}({\Bbb Z}))\otimes {\Bbb Q}_p$. Equation (8) implies that the 
restriction of this action of ${\Bbb T}(h)$ to $C^i(\Gamma,L_{\tilde{\lambda}}({\Bbb Z}))$ coincides with the action of ${\Bbb T}(h)$ on 
$C^i(\Gamma,L_{\tilde{\lambda}}({\Bbb Z}))$,  
i.e. ${\Bbb T}(h)$ on $C^i(\Gamma,L_{\tilde{\lambda}}({\Bbb Q}_p))$ is the linear extension of ${\Bbb T}(h)$ on $C^i(\Gamma,L_{\tilde{\lambda}}({\Bbb Z}))$. 
%
%
%
On the other hand, the definitions in (5.1.3) imply that
$$
L_{\tilde{\lambda}}(\frac{{\Bbb Z}_p}{(p^r)})=L_{\tilde{\lambda}}(\frac{{\Bbb Z}}{(p^r)})\otimes {\Bbb Z}_p, 
\quad L_{\tilde{\lambda}}({\Bbb Z}_p,r)=L_{\tilde{\lambda}}({\Bbb Z},r)\otimes {\Bbb Z}_p,
\quad L_{\tilde{\lambda}}(\frac{{\Bbb Z}_p}{(p^r)},r)=L_{\tilde{\lambda}}(\frac{{\Bbb Z}}{(p^r)},r)\otimes{\Bbb Z}_p,
$$
hence, we obtain

\begin{itemize}

\item $C^\bullet(\Gamma,L_{\tilde{\lambda}}({\Bbb Z}_p))=C^\bullet(\Gamma,L_{\tilde{\lambda}}({\Bbb Z}))\otimes{\Bbb Z}_p$

\item $C^\bullet(\Gamma,L_{\tilde{\lambda}}(\frac{{\Bbb Z}_p}{(p^r)}))=C^\bullet(\Gamma,L_{\tilde{\lambda}}(\frac{{\Bbb Z}}{(p^r)}))\otimes{\Bbb Z}_p$

\item $C^\bullet(\Gamma,L_{\tilde{\lambda}}({\Bbb Z}_p),r)=C^\bullet(\Gamma,L_{\tilde{\lambda}}({\Bbb Z},r))\otimes{\Bbb Z}_p$

\item $C^\bullet(L_{\tilde{\lambda}}(\frac{{\Bbb Z}_p}{(p^r)},r))=C^\bullet(L_{\tilde{\lambda}}(\frac{{\Bbb Z}}{(p^r)},r))\otimes{\Bbb Z}_p$

\item $C^\bullet(\Gamma, {\bf L}^{[r]}_{\tilde{\lambda}}({\Bbb Z}_p))
=C^\bullet(\Gamma,{\bf L}^{[r]}_{\tilde{\lambda}}({\Bbb Z}))\otimes{\Bbb Z}_p$.
%

\end{itemize}

The discussion in (5.7) shows that ${\Bbb T}(h)$ leaves all of the above complexes invariant and annihilates the second last one.  We obtain

\bigskip

{\bf Proposition. } {\it The normalized Hecke operator ${\Bbb T}(h)$ acts on the cohomology  groups

\medskip

$H^i(\Gamma,L_{\tilde{\lambda}}({\Bbb Z}_p))$, $H^i(\Gamma,L_{\tilde{\lambda}}(\frac{{\Bbb Z}_p}{(p^r)}))$, $H^i(\Gamma,L_{\tilde{\lambda}}({\Bbb Z}_p,r))$, 
$H^i(\Gamma,L_{\tilde{\lambda}}(\frac{{\Bbb Z}_p}{(p^r)},r))$, $H^i(\Gamma,{\bf L}^{[r]}_{\tilde{\lambda}}({\Bbb Z}_p))$

\medskip

and ${\Bbb T}(h)$ annihilates the second last one.
}

\section{Boundedness of slope spaces}

{\bf (6.1)  Slope subspaces. } We fix a finite dimensional ${\Bbb Q}_p$-vector space $V$ on which the operator $T$ acts and we assume that $T$ leaves a 
${\Bbb Z}_p$-lattice $V_{{\Bbb Z}_p}$ in $V$ invariant.
We denote by 
$m(X)\in {\Bbb Q}_p[X]$ the characteristic polynomial of $T$ acting on $V$. Since $T$ leaves $V_{{\Bbb Z}_p}$ invariant we know that $m(X)\in {\Bbb Z}_p[X]$. 
For any $\alpha\in{\Bbb Q}_{\ge 0}$ we set
$$
m_\alpha(X)=\prod_{\mu,\;v_p(\mu)\not=\alpha} (X-\mu),
$$ 
where $\mu$ runs over all roots of $m(X)$ in a splitting field for $T$ (counted with their multiplicities), which have $p$-adic value not equal to $\alpha$. We note 
that $m_\alpha(X)$ is contained in ${\Bbb Z}_p[X]$ because the condition "$v_p(\mu)\not=\alpha$" is ${\rm Gal}(\bar{\Bbb Q}_p/{\Bbb Q}_p)$-invariant and the roots $\mu$ 
are integers. 

\medskip

Let $K/{\Bbb Q}_p$ be a finite extension field with valuation $v_p$ and let ${\cal O}$ be the integers in $K$. We set $V_K=V\otimes K$ and 
$V_{\cal O}=V_{{\Bbb Z}_p}\otimes{\cal O}$ (hence, $V_{{\Bbb Q_p}}=V$). For any $\alpha\in{\Bbb Q}_{\ge 0}$ and we define the slope $\alpha$ subspaces $V_K^\alpha$ of $V_K$ resp. 
$V_{\cal O}^\alpha$ of $V_{\cal O}$ as $V_K^\alpha=m_\alpha(T) V_K$ resp. as $V_{\cal O}^\alpha=V_K^\alpha\cap V_ {\cal O}$. We note that $V_{\cal O}^\alpha$ is a lattice in $V_K^\alpha$ 
(i.e. $V_{\cal O}^\alpha$ contains a ${\cal O}$-basis which is a $K$-basis of $V_K^\alpha$) 
and $T$ leaves $V_K^\alpha$ and $V_{\cal O}^\alpha$ invariant. 
For brevity, we write $V^\alpha=V_{{\Bbb Q}_p}^\alpha$. We note that $V_K^\alpha$ is defined over ${\Bbb Q}_p$. 

\medskip

In this article, instead of a single slope subspace $V^{\alpha}_K$, we will consider the sum of all slope spaces whose slope is smaller than a given slope $\beta$. Thus, 
for any $\beta\in{\Bbb Q}_{\ge 0}$ we set
$$
V_K^{\le \beta}=\bigoplus_{0\le \alpha\le \beta} V_K^\alpha.
$$
Notice that $V_K^\alpha\subseteq V_K^{\le\alpha}$, hence, any bound on the dimension of $V_K^{\le\alpha}$ also yields a bound on the dimension of $V_K^\alpha$. We set 
$V_{\cal O}^{\le\beta}=V_K^{\le\beta}\cap V_{\cal O}$ and for simplicity we write $V^{\le\beta}=V^{\le \beta}_{{\Bbb Q}_p}$. Again, $V_{\cal O}^{\le\beta}$ is a lattice in $V_K^{\le\beta}$
%
%
%
%
%
which is invariant under $T$ and $V_K^{\le\beta}$ is defined over ${\Bbb Q}_p$. We define
$$
d(\beta)={\rm dim}\,V^{\le\beta}={\rm dim}\,V_K^{\le\beta}={\rm dim}\,V_{\cal O}^{\le\beta}.
$$

\bigskip

{\it Remark. } There are inclusions
$$
\bigoplus_{0\le\alpha\le \beta} m_\alpha(T) V_{\cal O}\subseteq \bigoplus_{0\le\alpha\le\beta} V_{\cal O}^\alpha\subseteq V_{\cal O}^{\le\beta},
$$
which in general are strict, i.e. equality does not hold in general in the above inclusions. 

%
%
%
%
%
%
%

\bigskip

{\bf (6.2)  A general estimate on the rank of slope subspaces. } We keep the notations from section (6.1); in particular, $T$ is an endomorphism of $V$ which 
leaves the lattice $V_{{\Bbb Z}_p}$ invariant.

\bigskip

{\bf (6.2.1)  Lemma. }{\it Let $\beta\in{\Bbb Q}_{\ge 0}$. Let $0\le \alpha_1<\alpha_2<\cdots<\alpha_s\le \beta$ be the non-trivial slopes appearing in $V^{\le \beta}$ and 
set $d_i={\rm dim}\,V^{\alpha_i}$. Then, for any integer $r$ we have the following bound for the dimensions $d_i$:
$$
\sum_{i=1}^s d_i(r-\alpha_i)\,\le\,
v_p\left( \#\,\left(\frac{V^{\le \beta}_{{\Bbb Z}_p}\otimes {\Bbb Z}_p/(p^r)}{{\rm ker}\,T|_{V^{\le \beta}_{{\Bbb Z}_p}\otimes {\Bbb Z}_p/(p^r)}}\right)\right).
$$
}

\medskip

(Note that the cardinality appearing on the right hand side is finite and a power of $p$.)




\bigskip

{\it Proof. } We abbreviate $d=d(\beta)$. We choose a splitting field $K/{\Bbb Q}_p$ for $T$ acting on $V$ and denote by ${\cal O}$ the ring of integers of $K$. In 
particular, $T$ is split on $V_K^{\alpha_i}$ for all $i=1,\ldots,s$, hence, there is a basis ${\cal B}_i$ of $V_K^{\alpha_i}$ such that $T|_{V_K^{\alpha_i}}$ is represented by a triangular matrix
$$
{\cal D}_{{\cal B}_i}(T|_{V_K^{\alpha_i}})=\left(\begin{array}{ccc}
\mu_{i,1}&&*\\
&\ddots&\\
&&\mu_{i,d_i}\\
\end{array}
\right).
$$ 
Here, $\mu_{i,1},\ldots,\mu_{i,d_i}\in{\cal O}$ are the eigenvalues of $T$ acting on $V_K^{\alpha_i}$ counted with their respective multiplicities. 
Thus, $v_p(\mu_{i,j})=\alpha_i$ for all $j=1,\ldots,d_i$ and since ${\rm det}\,T|_{V_K^{\le \beta}}=\prod_i {\rm det}\,T|_{V_K^{\alpha_i}}=\prod_{i,j}\mu_{i,j}$ 
we obtain
$$
v_{p}({\rm det}\,T|_{V_K^{\le \beta}})=\sum_{i=1,\ldots,s}\sum_{j=1,\ldots,d_i} v_{p}(\mu_{i,j})=\sum_{i=1,\ldots,s} d_i \alpha_i.\leqno(1)
$$
On the other hand, we choose a basis ${\cal C}$ of the lattice $V_{{\Bbb Z}_p}^{\le \beta}$ in $V^{\le\beta}$. Thus, ${\cal C}$ is a basis of $V^{\le\beta}$ and the 
representing matrix ${\cal D}_{\cal C}(T|_{V^{\le \beta}})$ has coefficients in ${\Bbb Z}_p$. Hence, the Theorem on elementary 
divisors implies that there are matrices $A,B\in{\rm GL}_d({\Bbb Z}_p)$ such that 
$$
A{\cal D}_{\cal C}(T|_{V^{\le \beta}})B=\left(\begin{array}{cccc}
\lambda_1&&\\
&\ddots&\\
&&\lambda_d\\
\end{array}
\right)\leqno(2)
$$ 
is diagonal ($\lambda_i\in{\Bbb Z}_p$). Equation (2) implies that $v_p({\rm det}\,T|_{V^{\le \beta}})=\sum_i v_p(\lambda_i)$, hence, using equation (1) we obtain
$$
\sum_{i=1}^d v_{p}(\lambda_i)=\sum_{i=1}^s d_i\alpha_i.\leqno(3)
$$
There is a basis $\{c_i\}_i$ and a basis $\{d_i\}_i$ of $V_{{\Bbb Z}_p}^{\le\beta}$ such that equation (2) gives the representing matrix of $T|_{V^{\le \beta}}$ with respect 
to the pair of basis $\{c_i\}_i$ and $\{d_i\}_i$:
%
%
%
$$
{\cal D}_{\{c_i\}_i,\{d_i\}_i}(T|_{V^{\le \beta}_{{\Bbb Z}_p}})=\left(\begin{array}{cccc}
\lambda_1&&\\
&\ddots&\\
&&\lambda_d\\
\end{array}
\right).
$$
Hence, we obtain
$$
\#\left(\frac{V_{{\Bbb Z}_p}^{\le\beta}\otimes{\Bbb Z}_p/(p^r)}{{\rm ker}\,T|_{V_{{\Bbb Z}_p}^{\le\beta}\otimes{\Bbb Z}_p/(p^r)}}\right)
=\#\,{\rm im}\,T|_{V_{{\Bbb Z}_p}^{\le\beta}\otimes{\Bbb Z}_p/(p^r)}
=\prod_{i=1}^d p^{{\rm max}\,(r-v_p(\lambda_i),0)}.
$$
Using equation (3) we obtain $\sum_{i=1}^s d_i(r-\alpha_i)=dr-\sum_{i=1}^s d_i\alpha_i=dr-\sum_{i=1}^d v_p(\lambda_i)=\sum_{i=1}^d r-v_p(\lambda_i)$, hence,
$$
p^{\sum_i d_i(r-\alpha_i)}=p^{\sum_i r-v_p(\lambda_i)}\,\big|\,p^{\sum_i {\rm max}\,(r-v_p(\lambda_i),0)}
=\#\left(\frac{V_{{\Bbb Z}_p}^{\le\beta}\otimes{\Bbb Z}_p/(p^r)}{{\rm ker}\,T|_{V_{{\Bbb Z}_p}^{\le\beta}\otimes{\Bbb Z}_p/(p^r)}}\right),
$$
which is the claim. (We note that the product $\alpha\,{\rm dim}\,V^\alpha\in{\Bbb N}$ is an integer (cf. [G-M], p. 797) and that $v_p(\lambda_i)\in{\Bbb Z}$ because 
$\lambda_i\in{\Bbb Z}_p$; hence, all exponents appearing in the above equation are integers.) 
Thus, the proof of the Lemma is complete.

\bigskip

{\bf (6.2.2) Corollary. } {\it Let $\beta\in{\Bbb Q}_{\ge 0}$ and let $r$ be any integer satisfying $r>\beta+1$. Then
$$
p^{d(\beta)}\, \Big | \, \#\,\left(\frac{V^{\le \beta}_{{\Bbb Z}_p}\otimes {\Bbb Z}_p/(p^r)}{{\rm ker}\,T|_{V^{\le \beta}_{{\Bbb Z}_p}\otimes {\Bbb Z}_p/(p^r)}}\right).
$$

}

{\it Proof. } Since $r>\beta+1\ge \alpha_i+1$ for all $i=1,\ldots,s$ we know that $\sum_{i=1}^s d_i(r-\alpha_i)\ge \sum_i d_i=d(\beta)$. Thus, the claim follows 
from the Lemma.

\bigskip

\bigskip

{\bf (6.3) Slope subspaces of the cohomology of arithmetic groups. } We turn to the situation of interest to us. We fix a reductive, ${\Bbb Q}$-split 
algebraic group $\tilde{\bf G}$, and we use the notations introduced in section (5.1). Thus, $\tilde{\bf G}(\bar{\Bbb Q})={\bf G}(\bar{\Bbb Q})\times_{\sf g}{\Bbb G}_m^a(\bar{\Bbb Q})$, where ${\bf G}={\bf G}_{\lambda_0}$ is the derived group of $\tilde{\bf G}$ and ${\sf g}\le {\bf G}(\bar{\Bbb Q})\cap {\Bbb G}_m^a(\bar{\Bbb Q})$ is a finite subgroup; moreover, $\tilde{\bf T}$ is a maximal split torus in $\tilde{\bf G}$ as defined in section (5.1.1), 
i.e. $\tilde{\bf T}(\bar{\Bbb Q})={\bf T}(\bar{\Bbb Q})\times_{\sf g}{\Bbb G}_m^a(\bar{\Bbb Q})$, where ${\bf T}$ is a maximal split torus in ${\bf G}$ as 
defined in (4.3). We fix the natural ${\Bbb Z}$-structure on ${\bf G}_{\lambda_0}$ (cf. (4.1)) and we select an arithmetic subgroup $\Gamma$ of ${\bf G}({\Bbb Z})$, 
which satisfies the following local condition at $p$:
$$
\Gamma\le K_*(p).
$$ 
(cf. section (4.2) for the definition of $K_*(p)$). Moreover, $\tilde{\lambda}=\lambda^\circ\otimes\kappa\in X(\tilde{\bf T})$ denotes any dominant weight as in equation (3) in section (5.1.2). As in section (5.2) 
we fix strictly positive integers $e_\alpha\in{\Bbb N}$ and an element $h\in \tilde{\bf T}({\Bbb Q})$ satisfying $v_p(\alpha^\circ(h))=e_\alpha$ for 
all $\alpha\in \Delta$ and we define the normalized Hecke operator 
$$
{\Bbb T}={\Bbb T}(h)={\tilde{\lambda}(h)} T(h).\leqno(4)
$$ 
We recall that ${\Bbb T}$ defines an operator on $H^i(\Gamma,L_{\tilde{\lambda}}({\Bbb Z}_p))$ and on $H^i(\Gamma,L_{\tilde{\lambda}}({\Bbb Z}_p)/L_{\tilde{\lambda}}({\Bbb Z}_p,r))$ (cf. (5.8)). We denote by $H^i(\Gamma,L_{\tilde{\lambda}}({\Bbb Z}_p))_{\rm int}$ the image of the canonical map
$$
H^i(\Gamma,L_{\tilde{\lambda}}({\Bbb Z}_p))\rightarrow H^i(\Gamma,L_{\tilde{\lambda}}({\Bbb Q}_p)).
$$
$H^i(\Gamma,L_{\tilde{\lambda}}({\Bbb Z}_p))_{\rm int}$ is a lattice in $H^i(\Gamma,L_{\tilde{\lambda}}({\Bbb Q}_p))$ and since the canonical map is 
${\Bbb T}$-equivariant, $H^i(\Gamma,L_{\tilde{\lambda}}({\Bbb Z}))_{\rm int}$ is a ${\Bbb T}$-stable lattice.
Thus, we are in the situation of (6.1) (with $V$ replaced by $H^i(\Gamma,L_{\tilde{\lambda}}({\Bbb Q}_p))$ and $V_{{\Bbb Z}_p}$ replaced by 
$H^i(\Gamma,L_{\tilde{\lambda}}({\Bbb Z}_p))_{\rm int}$) and following (6.1) we define the following slope subspaces:

\begin{itemize}

\item $H^i(\Gamma,L_{\tilde{\lambda}}({\Bbb Q}_p))^\alpha$ is the slope $\alpha$ subspace of $H^i(\Gamma,L_{\tilde{\lambda}}({\Bbb Q}_p))$ with respect to the 
normalized Hecke operator ${\Bbb T}$ 

\item $H^i(\Gamma,L_{\tilde{\lambda}}({\Bbb Z}_p))_{\rm int}^\alpha=H^i(\Gamma,L_{\tilde{\lambda}}({\Bbb Q}_p))^\alpha\cap H^i(\Gamma,L_{\tilde{\lambda}}({\Bbb Z}_p))_{\rm int}$


\item $H^i(\Gamma,L_{\tilde{\lambda}}({\Bbb Q}_p))^{\le\beta}=\bigoplus_{0\le\alpha\le\beta} H^i(\Gamma,L_{\tilde{\lambda}}({\Bbb Q}_p))^\alpha$

\item $H^i(\Gamma,L_{\tilde{\lambda}}({\Bbb Z}_p))_{\rm int}^{\le\beta}
=H^i(\Gamma,L_{\tilde{\lambda}}({\Bbb Q}_p))^{\le\beta}\cap H^i(\Gamma,L_{\tilde{\lambda}}({\Bbb Z}_p))_{\rm int}$.
\end{itemize}

We recall from (6.1) that $H^i(\Gamma,L_{\tilde{\lambda}}({\Bbb Z}_p))_{\rm int}^{\le\beta}$ is a lattice in $H^i(\Gamma,L_{\tilde{\lambda}}({\Bbb Q}_p))^{\le\beta}$.
In this final section we want to show that the dimensions
$$
d(\tilde{\lambda},i,\beta)={\rm dim}_{{\Bbb Q}_p}\,H^i(\Gamma,L_{\tilde{\lambda}}({\Bbb Q}_p))^{\le\beta} ={\rm dim}_{{\Bbb Z}_p}\,H^i(\Gamma,L_{\tilde{\lambda}}({\Bbb Z}_p))_{\rm int}^{\le\beta}
$$
are bounded by some constant $C(\beta,i,\Gamma,p)$ which does not depend on $\tilde{\lambda}$. 

\bigskip

{\bf (6.4) Proposition. } {\it There are ${\Bbb Z}_p/(p^r)$-submodules ${\cal X}\le H^i(\Gamma, L_{\tilde{\lambda}}({\Bbb Z}_p))_{\rm int}\otimes{\Bbb Z}_p/(p^r)$ and 
${\cal Y}\le H^i(\Gamma, {\bf L}^{[r]}_{\tilde{\lambda}}({\Bbb Z}_p))$ such that

- there is an exact sequence
$$
{\cal X}\stackrel{"\subseteq"}{\rightarrow} H^i(\Gamma, L_{\tilde{\lambda}}({\Bbb Z}_p))_{\rm int}\otimes{\Bbb Z}_p/(p^r)\rightarrow \frac{H^i(\Gamma, {\bf L}^{[r]}_{\tilde{\lambda}}({\Bbb Z}_p))}{{\cal Y}},
$$

- ${\cal X}$ and ${\cal Y}$ are invariant under ${\Bbb T}$ and ${\cal X}$ is annihilated by ${\Bbb T}$,

- the maps in the above sequence commute with the action of ${\Bbb T}$.
}

\medskip


{\it Proof. } For shortness, we set $H^i(V)=H^i(\Gamma,V)$ for any $\Gamma$-module $V$. We start from the short exact sequence
$$
0\rightarrow L_{\tilde{\lambda}}(\frac{{\Bbb Z}_p}{(p^r)},r)\stackrel{j}{\rightarrow} L_{\tilde{\lambda}}({\Bbb Z}_p)\otimes\frac{{\Bbb Z}_p}{(p^r)}\stackrel{\pi}{\rightarrow} 
\frac{L_{\tilde{\lambda}}({\Bbb Z}_p)}{L_{\tilde{\lambda}}({\Bbb Z}_p,r)}\rightarrow 0,\leqno(5)
$$
where $j$ is the inclusion and $\pi$ is defined by $\pi(v\otimes 1)=v+L_{\tilde{\lambda}}({\Bbb Z},r)$. 
%
%
From equation (5) we obtain a long exact cohomology sequence
$$
\ldots\rightarrow H^i(L_{\tilde{\lambda}}(\frac{{\Bbb Z}_p}{(p^r)},r))\stackrel{H(j)}{\rightarrow} H^i(L_{\tilde{\lambda}}({\Bbb Z}_p)\otimes\frac{{\Bbb Z}_p}{(p^r)})
\stackrel{H(\pi)}{\rightarrow} H^i({\bf L}^{[r]}_{\tilde{\lambda}}({\Bbb Z}_p))\rightarrow \cdots.\leqno(6)
$$
We note that ${\Bbb T}$ acts on all cohomology groups appearing in the above sequence (cf. (5.8)). We decompose $\Gamma h \Gamma=\bigcup_i \Gamma h\gamma_i$. 
(5.4) Lemma implies (after tensoring with ${\Bbb Z}_p$) that 
$\tilde{\lambda}(h)\pi_{\tilde{\lambda}}(h\gamma_i)^{-1}$ leaves $L_{\tilde{\lambda}}({\Bbb Z}_p)$ and $L_{\tilde{\lambda}}({\Bbb Z}_p,r)$ invariant. Hence, 
it also leaves the $\pmod{p^r}$-reduction $L_{\tilde{\lambda}}({\Bbb Z}_p/(p^r))$ of $L_{\tilde{\lambda}}({\Bbb Z}_p)$ and the image 
$L_{\tilde{\lambda}}({\Bbb Z}_p/(p^r),r)$ of $L_{\tilde{\lambda}}({\Bbb Z}_p,r)$ in $L_{\tilde{\lambda }}({\Bbb Z}_p)\otimes {\Bbb Z}_p/(p^r)$ under the map $v\mapsto v\otimes 1$ invariant. 
Thus, it acts on 
all components of the exact sequence in equation (5) and since it is obvious that 
the mappings $j$ and $\pi$ in equation (5) commute with $\tilde{\lambda}(h)\pi_{\tilde{\lambda}}(h\gamma_i)^{-1}$ it follows that the mappings $H(j)$, $H(\pi)$ 
in equation (6) commute with the action of the Hecke operator ${\Bbb T}$ (cf. [K-P-S], Theorem 1.3.1, p. 227; note that we normalized the action of the 
Hecke operator by multiplying with $\tilde{\lambda}(h)$). The short exact sequence
$$
0\rightarrow L_{\tilde{\lambda}}({\Bbb Z}_p)\stackrel{\cdot p^r}{\rightarrow}L_{\tilde{\lambda}}({\Bbb Z}_p)\stackrel{\rm pr}{\rightarrow} L_{\tilde{\lambda}}({\Bbb Z}_p)/p^r L_{\tilde{\lambda}}({\Bbb Z}_p)\rightarrow 0,
$$
where ${\rm pr}$ is the natural projection, yields a long exact sequence
$$
\cdots\rightarrow H^i(L_{\tilde{\lambda}}({\Bbb Z}_p))\stackrel{\cdot p^r}{\rightarrow}H^i(L_{\tilde{\lambda}}({\Bbb Z}_p))\stackrel{H(\rm pr)}{\rightarrow} 
H^i(L_{\tilde{\lambda}}({\Bbb Z}_p)/p^r L_{\tilde{\lambda}}({\Bbb Z}_p))\rightarrow\cdots
$$
hence, we obtain an embedding
$$
H^i(L_{\tilde{\lambda}}({\Bbb Z}_p))\otimes\frac{{\Bbb Z}_p}{(p^r)}\stackrel{H(\rm pr)}{\hookrightarrow} 
H^i(L_{\tilde{\lambda}}({\Bbb Z}_p)\otimes\frac{{\Bbb Z}_p}{(p^r)}).\leqno(7)
$$
We identify $H^i(L_{\tilde{\lambda}}({\Bbb Z}_p))\otimes{\Bbb Z}_p/(p^r)$ with its image under $H({\rm pr})$. Restricting the morphism $H(\pi)$ to $H^i(L_{\tilde{\lambda}}({\Bbb Z}_p))\otimes{\Bbb Z}_p/(p^r)$ we obtain from equation (6) the exact sequence
$$
X\stackrel{H(j)}{\rightarrow} H^i(L_{\tilde{\lambda}}({\Bbb Z}_p))\otimes\frac{{\Bbb Z}_p}{(p^r)}
\stackrel{H(\pi)}{\rightarrow} H^i({\bf L}^{[r]}_{\tilde{\lambda}}({\Bbb Z}_p)),\leqno(8)
$$
where
$$
X=H^i(L_{\tilde{\lambda}}(\frac{{\Bbb Z}_p}{(p^r)},r))\cap H(j)^{-1}\big( H^i(L_{\tilde{\lambda}}({\Bbb Z}_p))\otimes\frac{{\Bbb Z}_p}{(p^r)}\big).
$$
We denote by ${\cal T}_\lambda$ the torsion submodule of $H^i(L_{\tilde{\lambda}}({\Bbb Z}_p))$ and we set ${\cal T}_{\lambda,r}={\cal T}_\lambda\otimes {\Bbb Z}_p/(p^r)$.
From equation (8) we further obtain the exact sequence
$$
\frac{X}{H(j)^{-1}({\cal T}_{\lambda,r})}\stackrel{H(j)}{\rightarrow} \frac{H^i(L_{\tilde{\lambda}}({\Bbb Z}_p))\otimes{\Bbb Z}_p/(p^r)}{{\cal T}_{\lambda,r}}
\stackrel{H(\pi)}{\rightarrow} \frac{H^i({\bf L}^{[r]}_{\tilde{\lambda}}({\Bbb Z}_p))}{H(\pi)({\cal T}_{\lambda,r})}.\leqno(9)
$$
We note that ${\Bbb T}$ acts on all modules appearing in equation (9), because ${\Bbb T}$ leaves the torsion ${\cal T}_\lambda$ invariant. We 
set ${\cal X}=X/H(j)^{-1}({\cal T}_{\lambda,r})$. Since ${\Bbb T}$ annihilates $H^i(L_{\tilde{\lambda}}(\frac{{\Bbb Z}_p}{(p^r)},r))$ (cf. (5.8)) 
we see that ${\Bbb T}$ annihilates $X$ and, hence, ${\cal X}$. Since moreover, 
$$
\frac{H^i(L_{\tilde{\lambda}}({\Bbb Z}_p))\otimes{\Bbb Z}_p/(p^r)}{{\cal T}_{\lambda,r}}\cong H^i(L_{\tilde{\lambda}}({\Bbb Z}_p))_{\rm int}\otimes{\Bbb Z}_p/(p^r)
$$ 
we see that equation (9) yields the claim after replacing ${\cal X}$ by its image under $H(j)$; note that ${\Bbb T}$ also annihilates 
the image of ${\cal X}$ under $H(j)$ because $H(j)$ and ${\Bbb T}$ commute. Thus, the Proposition is proven.

\bigskip

{\bf (6.5) Proposition. }{\it Let $\beta\in{\Bbb Q}_{\ge 0}$ and choose an integer $r$ bigger than $\beta+1$. Then, for all dominant weights $\tilde{\lambda}\in X(\tilde{\bf T})$ we have
$$
p^{d(\tilde{\lambda},i,\beta)}\,\Big|  \, \#\,(H^i(\Gamma,{\bf L}^{[r]}_{{\lambda}}({\Bbb Z}_p))).
$$

}

\medskip

{\it Proof. } We let $\tilde{\lambda}\in X(\tilde{\bf T})$ be an arbitrary dominant weight, i.e. 
$\tilde{\lambda}=\lambda^\circ\otimes\kappa$ where $\lambda^\circ$ corresponds to a dominant weight $\lambda=\sum_{\alpha\in\Delta} m_\alpha\omega_\alpha$, $m_\alpha\in{\Bbb N}_0$, in 
$\Gamma_{\lambda_0}$ (cf. section (4.3) equation (5)). Since $H^i(\Gamma,L_{\tilde{\lambda}}({\Bbb Z}_p))^{\le\beta}_{\rm int}$ is saturated in $H^i(\Gamma,L_{\tilde{\lambda}}({\Bbb Z}_p))_{\rm int}$
the canonical map
$$
H^i(\Gamma,L_{\tilde{\lambda}}({\Bbb Z}_p))_{\rm int}^{\le\beta}\otimes\frac{{\Bbb Z}_p}{(p^r)}\rightarrow H^i(\Gamma,L_{\tilde{\lambda}}({\Bbb Z}_p))_{\rm int}\otimes\frac{{\Bbb Z}_p}{(p^r)}
$$
is injective.
%
%
Thus, we may restrict the exact sequence in (6.4) Proposition to 
the submodule $H^i(\Gamma,L_{\tilde{\lambda}}({\Bbb Z}_p))^{\le\beta}_{\rm int}\otimes{\Bbb Z}_p/(p^r)$ of $H^i(\Gamma,L_{\tilde{\lambda}}({\Bbb Z}_p))_{\rm int}\otimes{\Bbb Z}_p/(p^r)$. More 
precisely, we set
$$
{\cal X}'={\cal X}\cap H^i(\Gamma,L_{\tilde{\lambda}}({\Bbb Z}_p))_{\rm int}^{\le \beta}\otimes{\Bbb Z}_p/(p^r)
$$ 
and obtain an exact sequence
$$
{\cal X}'\stackrel{"\subseteq"}{\rightarrow} H^i(\Gamma,L_{\tilde{\lambda}}({\Bbb Z}_p))^{\le\beta}_{\rm int}\otimes{\Bbb Z}_p/(p^r)
\rightarrow \frac{H^i(\Gamma,{\bf L}^{[r]}_{\tilde{\lambda}}({\Bbb Z}_p))}{\cal Y}.
$$
Since ${\Bbb T}$ annihilates the submodule ${\cal X}$ of $ H^i(\Gamma,L_{\tilde{\lambda}}({\Bbb Z}_p))_{\rm int}\otimes{\Bbb Z}_p/(p^r)$ (cf. (6.4) Proposition) 
we know that ${\cal X}'$ is contained in the kernel of ${\Bbb T}|_{H^i(\Gamma,L_{\tilde{\lambda}}({\Bbb Z}_p))_{\rm int}^{\le\beta}\otimes{\Bbb Z}_p/(p^r)}$. Hence, we obtain
$$
\#\,\left(\frac{H^i(\Gamma,L_{\tilde{\lambda}}({\Bbb Z}_p))_{\rm int}^{\le\beta}\otimes{\Bbb Z}_p/(p^r)}{{\rm ker}\,{\Bbb T}|_{H^i(\Gamma,L_{\tilde{\lambda}}({\Bbb Z}_p))_{\rm int}^{\le\beta}\otimes{\Bbb Z}_p/(p^r)}}\right)   \,  \Big| \,
\# \, (H^i(\Gamma,{\bf L}^{[r]}_{\tilde{\lambda}}({\Bbb Z}_p)))
$$
(note that all (cohomology) groups appearing in the above equation are finite groups). 
(6.2.2) Corollary then implies that 
$$
p^{d(\tilde{\lambda},i,\beta)}\,\Big|  \, \#\,(H^i(\Gamma,{\bf L}^{[r]}_{\tilde{\lambda}}({\Bbb Z}_p))).
$$
Since $\Gamma\le{\bf G}({\Bbb Z})$ we know that ${\bf L}^{[r]}_{\tilde{\lambda}}({\Bbb Z}_p)={\bf L}^{[r]}_\lambda({\Bbb Z}_p)$ (cf. (5.1.3))
and we obtain
$$
p^{d(\tilde{\lambda},i,\beta)}\,\Big|  \, \#\,(H^i(\Gamma,{\bf L}^{[r]}_{{\lambda}}({\Bbb Z}_p))).
$$
Thus, the proof of the Proposition is complete.

\bigskip

{\bf (6.6) Boundedness of the slope subspaces. } (6.5) Proposition and (4.4) Proposition together yield our main result.

\bigskip

{\bf Theorem. }{\it Let $\beta\in{\Bbb Q}_{\ge 0}$ and choose $i\in{\Bbb N}_0$. There is $C=C(\beta,i,\Gamma,p)\in{\Bbb N}$ only depending on 
$\beta,i,\Gamma,p$ such that  
$$
{\rm dim}\,H^i(\Gamma,L_{\tilde{\lambda}}({\Bbb Q}_p))^{\le\beta}\le C.
$$
for all dominant weights ${\tilde{\lambda}}\in X(\tilde{\bf T})$.
}

\medskip

{\it Proof. } We set ${\cal F}=\Gamma_{\rm sc}/\Gamma_{\lambda_0}$, i.e. ${\cal F}$ is a quotient of the fundamental group of ${\bf G}$, and we let
$$
{\sf f}=\#({\cal F}).
$$
Thus, ${\sf f}\omega_\alpha\in \Gamma_{\lambda_0}$ for all fundamental weights $\omega_\alpha$. We fix an integer $r$ bigger than $\beta+1$ and we let 
$\tilde{\lambda}=\lambda^\circ\otimes\kappa\in X(\tilde{\bf T})$ be an arbitrary dominant weight; thus, the weight $\lambda^\circ$ corresponds to a  
dominant weight $\lambda=\sum_{\alpha\in\Delta}m_\alpha\omega_\alpha$ in $\Gamma_{\lambda_0}$. We denote by $\Delta(\lambda,r)$ the set of all simple roots such that 
twice the distance of $\lambda$ to the corresponding wall of the Weyl chamber is larger than $r$, i.e.
$$
\Delta(\lambda,r)=\{\alpha\in\Delta: m_\alpha=\lambda(h_\alpha)>r\}.
$$
Substracting from $\lambda$ suitable integral multiples of the weights $p^{\lceil\frac{p}{p-1}r \rceil}{\sf f}\omega_\alpha\in \Gamma_{\lambda_0}$ for $\alpha\in\Delta(\lambda,r)$, we obtain 
an integral and dominant weight $\lambda'=\sum_{\alpha\in\Delta} m_\alpha'\omega_\alpha\in\Gamma_{\lambda_0}$ ($m_\alpha'\in{\Bbb N}_0$) such that

\begin{itemize}

\item $m_\alpha'=m_\alpha$ for all $\alpha\in \Delta-\Delta(\lambda,r)$

\item $r<m_\alpha'\le r+p^{\lceil\frac{p}{p-1}r \rceil}{\sf f}$ for all $\alpha\in\Delta(\lambda,r)$

\item $\lambda'\equiv \lambda\pmod{p^{\lceil\frac{p}{p-1}r \rceil} \Gamma_{\rm sc}}$ \quad (note that $\omega_\alpha\in\Gamma_{\rm sc}$).

\end{itemize}

In particular, $m_\alpha,\,m_\alpha'>r$ for all $\alpha\in\Delta(\lambda,r)$ and (4.4) Proposition (with ${\cal T}=\Delta(\lambda,r)$) implies that 
${\bf L}^{[r]}_\lambda({\Bbb Z}_p)\cong {\bf L}^{[r]}_{\lambda'}({\Bbb Z}_p)$ as $\Gamma$-modules. 
Hence, we obtain 
$$
\#(H^i(\Gamma,{\bf L}^{[r]}_\lambda({\Bbb Z}_p)))= \#(H^i(\Gamma,{\bf L}^{[r]}_{\lambda'}({\Bbb Z}_p))).\leqno(10) 
$$
We define $\Lambda(r)\subseteq \Gamma_{\rm sc}$ as the set consisting of all integral weights $\mu=\sum_{\alpha\in\Delta} n_\alpha \omega_\alpha$ satisfying 
$0\le n_\alpha\le r+p^{\lceil\frac{p}{p-1}r \rceil}{\sf f}$ for all $\alpha\in \Delta$. 
Hence, $\lambda'$ is contained in $\Lambda(r)$ and equation (10) therefore yields
$$
\#(H^i(\Gamma,{\bf L}^{[r]}_\lambda({\Bbb Z}_p)))
\le C:={\rm max}_{\lambda'\in\Lambda(r)} \#(H^i(\Gamma,{\bf L}^{[r]}_{\lambda'}({\Bbb Z}_p))).
$$
The finiteness of $L_{\lambda'}({\Bbb Z}_p)/L_{\lambda'}({\Bbb Z}_p,r)$ implies that the cohomology $H^i(\Gamma,L_{\lambda'}({\Bbb Z}_p)/L_{\lambda'}({\Bbb Z}_p,r))$ is an 
abelian group of finite cardinality (note that we assume $\Gamma$ to be arithmetic). 
Since $\Lambda(r)$ is a finite set we obtain that $C$ as defined above is a finite constant which does not depend on 
$\tilde{\lambda}$ (it only depends on ${\bf G}$, $p$, $\Gamma$, $\beta$ and $i$). (6.5) Proposition now yields
$$
p^{d(\tilde{\lambda},i,\beta)}\le C.
$$
Since $\tilde{\lambda}$ was arbitrary, this is the claim and the proof of the theorem is complete.

\bigskip

\bigskip

{\Large\bf References}

[B] Borel, A. et al., Seminar on Algebraic groups and related finite groups, LNM {\bf 131}, Springer 1970

[Bu] Buzzard, K., Families of modular forms, J. de Th{\'e}orie de Nombres de Bordeaux {\bf 13} (2001) 43 - 52

[G-M] Gouv{\^e}a, F., Mazur, B., Families of modular eigenforms, Math. of Comp. {\bf 58}(198) (1992) 793 - 805

[H] Harder, G., Interpolating coefficient systems and $p$-ordinary cohomology of arithmetic groups, preprint 2010


[Ha] Hall, B., Lie groups, Lie algebras and representations, Springer 2003

[Hi] Hida, H., Elementary Theory of $L$-functions and Eisenstein series, London Mathematical Society Student Texts {\bf 26}, Cambridge University Press 1993

[Hu] Humphreys, J., Introduction to Lie algebras and representation theory, Springer 1972

[K-P-S], Kuga, M., Parry, W., Sah, C., Group cohomology and Hecke operators, in: Manifolds and Lie groups, ed. Hano, J. et
al, p. 223 - 266, Birkhauser 1981

[M 1] Mahnkopf, J., On $p$-adic families of modular forms, Proc. of RIMS Conf. on "Automorphic Forms, Automorphic representations and Related Toopics", 2010, pp. 93 - 108

[M 2] -, $P$-adic families of modular forms, to appear in: Oberwolfach reports 2011

[P] Pande, A., Local constancy of dimensions of Hecke eigenspaces of automorphic forms, Journal of Number Theory {\bf 129} (2009)  15-27 

[S] Springer, T., Linear algebraic groups, 2nd ed., 2nd printing, Birkh{\"a}user, Boston 2008

\newpage

\end{document}